\documentclass[11pt,reqno]{amsart}

\usepackage{amssymb}
\usepackage{mathdots} %
\usepackage{cancel} %
\usepackage[shortlabels]{enumitem} %
\usepackage{geometry}
\usepackage{verbatim} %

\usepackage{hyperref}
\usepackage{xcolor}
\definecolor{darkgreen}{rgb}{0,0.45,0}
\hypersetup{colorlinks,urlcolor=blue,citecolor=darkgreen,linkcolor=blue,linktocpage=false}
\usepackage[abbrev]{amsrefs}
\usepackage[capitalise,noabbrev]{cleveref}
\usepackage[all]{xy}
\newdir{ >}{{}*!/-8pt/@{>}} %
\newcommand{\cxymatrix}[1]{\vcenter{\xymatrix{#1}}} %
\SilentMatrices

\usepackage{ifpdf}
\ifpdf
\else
\usepackage{breakurl}
\fi

\usepackage{mathtools} %

\numberwithin{equation}{section}

\theoremstyle{plain}
\newtheorem{thm}[equation]{Theorem}
\newtheorem*{thm*}{Theorem} %
\newtheorem{prop}[equation]{Proposition}
\newtheorem{lem}[equation]{Lemma}
\newtheorem{cor}[equation]{Corollary}
\newtheorem*{cor*}{Corollary} %

\crefname{prop}{Proposition}{Propositions}

\theoremstyle{definition}
\newtheorem{defn}[equation]{Definition}
\newtheorem{ex}[equation]{Example}

\crefname{ex}{Example}{Examples}
\crefname{lem}{Lemma}{Lemmas}

\theoremstyle{remark}

\newtheorem{rem}[equation]{Remark}

\makeatletter
\newcommand{\Eqref}[1]{\textup{\tagform@{\ref*{#1}}}}
\makeatother

\newcommand{\anti}{\ar@{}[dr]|{\;\;\;\boxed{\textup{\scriptsize -1}}}} %

\newcommand{\Def}[1]{\textbf{\boldmath{#1}}} %

\newcommand{\F}{\mathbb{F}}
\newcommand{\Z}{\mathbb{Z}}

\newcommand{\DZ}{D(\Z)}

\newcommand{\cat}[1]{\mathbf{\mathcal{#1}}} %

\newcommand{\bmat}[1]{\left[\!\!\!\begin{array}{#1}} %
\newcommand{\emat}{\end{array}\!\!\right]} %

\newcommand{\smat}[1]{\left[\begin{smallmatrix*}[r]#1\end{smallmatrix*}\right]} %

\newcommand{\cT}{\cat{T}}
\newcommand{\cTo}{\cat{T}^{\opp}}

\newcommand{\al}{\alpha}
\newcommand{\be}{\beta}
\newcommand{\de}{\delta}

\newcommand{\del}{\partial}
\newcommand{\ep}{\epsilon}
\newcommand{\ga}{\gamma}
\newcommand{\Ga}{\Gamma}

\newcommand{\phy}{\varphi}
\newcommand{\Si}{\Sigma}

\newcommand{\te}{\theta}

\newcommand{\circar}{\ar|(.48)*\cir<1.7pt>{}} %

\newcommand{\ral}{\xrightarrow} %
\renewcommand{\implies}{\Longrightarrow}
\renewcommand{\impliedby}{\Longleftarrow}
\renewcommand{\iff}{\Longleftrightarrow}

\newcommand{\op}{\oplus}

\newcommand{\x}{\times}

\newcommand{\hatt}{\widehat}
\newcommand{\lan}{\left\langle}
\newcommand{\ol}{\overline}
\newcommand{\ran}{\right\rangle}
\newcommand{\tild}{\widetilde}

\newcommand{\StMod}{\mathrm{StMod}}

\DeclareMathOperator{\Ext}{Ext}

\newcommand{\inc}{\mathrm{inc}}
\newcommand{\opp}{\mathrm{op}}
\newcommand{\proj}{\mathrm{proj}}

\def\noteson{\gdef\note##1{\noindent{\color{blue}[##1]}}}

\noteson

\begin{document}

\title{On good morphisms of exact triangles} 
\date{\today}

\author{J. Daniel Christensen}
\address{Department of Mathematics\\
University of Western Ontario\\
London, Ontario, N6A 5B7\\
Canada}             
\email{jdc@uwo.ca}

\author{Martin Frankland}             
\address{University of Regina\\
3737 Wascana Parkway\\
Regina, Saskatchewan, S4S 0A2\\
Canada}
\email{Martin.Frankland@uregina.ca}

\begin{abstract}
In a triangulated category, cofibre fill-ins always exist. Neeman showed that there is always at least one ``good'' fill-in, i.e., one whose mapping cone is exact. Verdier constructed a fill-in of a particular form in his proof of the $4 \x 4$ lemma, which we call ``Verdier good''. We show that for several classes of morphisms of exact triangles, the notions of good and Verdier good agree. 
We prove a lifting criterion for commutative squares in terms of (Verdier) good fill-ins. 
Using our results on good fill-ins,
we also prove a pasting lemma for homotopy cartesian squares.
\end{abstract}



\thanks{\copyright 2021. This manuscript version is made available under the CC-BY-NC-ND 4.0 license \url{http://creativecommons.org/licenses/by-nc-nd/4.0/} }

\maketitle

\tableofcontents

\section{Introduction}

In this paper, we study various classes of maps between triangles in a
triangulated category $\cT$, and use our results to prove lifting
criteria for commutative squares, a pasting lemma for homotopy
cartesian squares, and a variety of results comparing these classes
of maps.

In order to explain our results, we make a few definitions.
A \Def{map of triangles} is a commutative diagram
\[
  \cxymatrix{
    X \ar[r]^u \ar[d]_f & Y \ar[d]^g \ar[r]^{v} & Z \ar[d]^h \ar[r]^w & \Sigma X \ar[d]^{\Sigma f} \\
    X' \ar[r]_{u'} & Y' \ar[r]_{v'} & Z' \ar[r]_{w'} & \Sigma X' \\
  }
\]
in which the rows are (exact) triangles, i.e., those specified by the triangulation.
We often write $(f,g,h)$ for such a map.

Following \cite{Neeman91}, we say that a map $(f,g,h)$ of triangles is
\Def{middling good} if it can be extended to a $4 \x 4$ diagram as in
\cref{def:middling-good}.
It is \Def{good} if its mapping cone (\cref{def:MappingCone}) is exact.
Neeman showed that every good map of triangles is middling good, but
that there exist maps which are middling good but not good, as well as
maps which are not middling good.

If $f$ and $g$ are given, a \Def{fill-in} is a map $h$ making the
rest of the diagram commute.  The fill-in is \Def{(middling) good}
if the resulting map $(f,g,h)$ of triangles is.

In this paper, we introduce a new condition:
a map $(f,g,h)$ of triangles is \Def{Verdier good} if $h$ is
constructed as in Verdier's proof that middling good fill-ins exist;
see \cref{def:VerdierGood}.
As a result, every Verdier good morphism is middling good.
\cref{lem:VerdierHaynes} provides a more conceptual characterization of Verdier good maps
that was suggested to us by Haynes Miller.

A natural question is whether either of good or Verdier good implies the other.
While we do not answer this question in general, we show that it
holds in a number of special situations.
For example, for maps of triangles with at most one non-zero component,
we show that being good, Verdier good, middling good, and nullhomotopic
are equivalent, and also show that these are equivalent to a certain
Toda bracket containing zero (see \cref{pr:ZeroVerdierGood}).

The case where a map of triangles has at most two non-zero components
is more interesting, and leads to our new lifting criteria.
For maps of the form $(0,g,h)$ and $(f,0,h)$, we show that being
good, Verdier good and nullhomotopic are equivalent (\cref{pr:0ghVerdierGood}).
Using this, we prove the following lifting criteria:

\begin{thm*}[{\cref{cor:lifting}}]
Given a solid arrow commutative square
\[
  \cxymatrix{
    X \ar[r]^u \ar[d]_f & Y \ar[d]^g \ar@{-->}[dl]_k \\
    X' \ar[r]_{u'} & Y' , \\
  }
\]
choose extensions of the rows to cofibre sequences with cofibres $Z$ and $Z'$,
respectively. The following are equivalent.
\begin{enumerate}
\item There exists a lift $k \colon Y \to X'$.
\item The map $0 \colon Z \to Z'$ is a fill-in and the map $(f,g,0)$ of triangles is good.
\item The map $0 \colon Z \to Z'$ is a fill-in and the map $(0,f,g)$ of rotated triangles is Verdier good.
\end{enumerate}
\end{thm*}

These criteria can be used to define obstruction classes as subsets of $[Z, Z']$,
which contain zero if and only if the lift exists.
This is described after \cref{th:obstruction} and \cref{cor:lifting},
and is related to the work in~\cite{ChristensenDI04}.

The rotation is needed in last case of the theorem because \cref{pr:0ghVerdierGood}
made no claims about maps of the form $(f,g,0)$, and we do not know if
being Verdier good is invariant under rotation.
Under a mild assumption, %
we show in \cref{pr:fg0Good} %
that for maps of the form $(f,g,0)$, good implies Verdier good. 
This extra assumption holds in most triangulated categories that arise in nature,
e.g., in the homotopy category of any stable model category or
of any complete and cocomplete stable $\infty$-category; see \cref{ex:n-angulated}.

Continuing the theme mentioned earlier, we give other classes of maps
of triangles for which goodness and Verdier goodness agree.
\cref{pr:IncUpperTriangular} shows that this holds for certain
split inclusions or split projections.
\cref{pr:ContractibleVerdierGood} shows that when the domain or codomain triangle
is contractible, every map is both good and Verdier good.
Goodness is invariant under chain homotopy, but we do not know
whether the same is true for Verdier goodness.
In fact, we do not know whether every nullhomotopic map is Verdier good.

\cref{pr:FillInRotated} studies certain maps from a triangle to its
rotation, and again shows that goodness and Verdier goodness agree on this class.
We use this to prove \cref{cor:examples}, which says:
there exist maps of triangles which are not middling good;
there exist maps of triangles which are middling good but neither
good nor Verdier good;
and middling goodness is not invariant under chain homotopy.
We give explicit examples in the same section.

The final family of maps we consider are maps of the form $(1,g,h)$,
where we have that $X = X'$.
In \cref{prop:HoCartFillin}, we show that such a map is good if and only if
the center square is homotopy cartesian (\cref{def:HoCart}) with a specific differential, %
reformulating a result of Neeman \cite{Neeman01}. 
We also characterize when the center square is homotopy cartesian. %
The main result of this section, \cref{pr:Pasting},
states that the pasting of two homotopy cartesian squares is again
homotopy cartesian.
While the statement does not mention good morphisms, the proof makes
essential use of them.
We then use this to strengthen a result of \cite{Vaknin01} about
a natural candidate triangle associated to the pasting.
We also characterize which maps of the form $(1,g,h)$ are Verdier good, %
and show that in this case good implies Verdier good. %

\subsection*{Open questions}

\begin{enumerate}
\item Is Verdier goodness invariant under chain homotopy? In particular, is every nullhomotopic map Verdier good?
\item Is Verdier goodness invariant under rotation?
\item Do either of good or Verdier good imply the other? In particular, given $f$ and $g$, is there a fill-in $h$ that is simultaneously good and Verdier good?
\end{enumerate}

\subsection*{Organization}

We begin in \cref{se:good} by giving background on maps of triangles and defining
the notion of a good map of triangles.
In \cref{se:vgood}, we introduce the notion of a Verdier good map of
triangles, give a characterization due to Haynes Miller, and
prove some first results about Verdier good maps.
In \cref{se:obstructions}, we give our first lifting criterion,
in terms of good fill-ins.
The remaining sections are organized by the family of maps that they study.
In \cref{se:incproj}, we study various split inclusions and split
projections, and show that every map of triangles is a composite of
two maps that are good and Verdier good. 
We study maps of the form $(1,g,h)$ as well as homotopy cartesian squares
in \cref{se:hocart}.
This section contains our results about pasting such squares, and also has
examples illustrating the subtlety in the choice of differential.
In \cref{se:zero} we focus on maps that have at most two non-zero components. %
Our final section, \cref{se:rotations}, deals with certain maps from
a triangle to its rotation, which allows us to produce a number of
interesting examples.

\subsection*{Conventions}

We work in a fixed triangulated category $\cT$ throughout, with the self-equivalence denoted $\Si \colon \cT \ral{\simeq} \cT$. 
We suppress the natural isomorphisms $\Si \Si^{-1} \cong 1$ and $\Si^{-1} \Si \cong 1$ from the notation. %
A \Def{candidate triangle} is a diagram in $\cT$ of the form 
\[
\cxymatrix{
X \ar[r]^-{u} & Y \ar[r]^-{v} & Z \ar[r]^-{w} & \Si X \\
}
\]
satisfying $vu = 0$, $wv = 0$, and $(\Si u)w = 0$ \cite{Neeman01}*{Definition~1.1.1}. 
A candidate triangle that belongs to the class specified by the triangulated structure is often called a \emph{distinguished triangle}; we will say \Def{exact triangle}, or simply \Def{triangle} for short. 

\subsection*{Acknowledgments}

We thank Haynes Miller, Amnon Neeman, and Marius Thaule for helpful discussions. 
We also thank the referee for a careful reading of the paper.

We acknowledge the support of the Natural Sciences and Engineering Research Council of Canada (NSERC). Cette recherche a \'et\'e financ\'ee par le Conseil de recherches en sciences naturelles et en g\'enie du Canada (CRSNG). RGPIN-2016-04648 (Christensen), RGPIN-2019-06082 (Frankland).

\section{Good morphisms}\label{se:good}

In this section, we provide some background on the notion of a ``good morphism''
between triangles, which was first introduced by Neeman in~\cite{Neeman91}.

\begin{defn}
A \Def{map of candidate triangles} is a commutative diagram
\begin{equation}\label{eq:map-of-triangles}
  \cxymatrix{
    X \ar[r]^u \ar[d]_f & Y \ar[d]^g \ar[r]^{v} & Z \ar[d]^h \ar[r]^w & \Sigma X \ar[d]^{\Sigma f} \\
    X' \ar[r]_{u'} & Y' \ar[r]_{v'} & Z' \ar[r]_{w'} & \Sigma X' \\
  }
\end{equation}
in which the rows are candidate triangles.
We often write $(f,g,h)$ for such a map.
We will refer to~\eqref{eq:map-of-triangles} many times in this paper.
\end{defn}

\begin{defn}\label{def:middling-good}
A map of triangles as in~\eqref{eq:map-of-triangles} is \Def{middling good}
if it extends to a $4 \x 4$ diagram
\begin{equation*}
\cxymatrix{
X \ar[d]_f \ar[r]^-{u} & Y \ar[d]^g \ar[r]^{v} & Z \ar[d]^h \ar[r]^-{w} & \Sigma X \ar[d]^{\Sigma f} \\
X' \ar[d]_{f'} \ar[r]^-{u'} & Y' \ar[d]^{g'} \ar[r]^-{v'} & Z' \ar[d]^{h'} \ar[r]^{w'} & \Sigma X' \ar[d]^{\Si f'} \\
X'' \ar[d]_{f''} \ar[r]^-{u''} & Y'' \ar[d]^{g''} \ar[r]^-{v''} & Z'' \anti \ar[d]^{h''} \ar[r]^-{w''} & \Sigma X'' \ar[d]^{\Sigma f''} \\
\Si X \ar[r]_-{\Si u} & \Si Y \ar[r]_-{\Si v} & \Si Z \ar[r]_-{\Si w} & \Sigma^2 X \\
}
\end{equation*}
where the first three rows and columns are exact, and the bottom-right square anticommutes as indicated.
\end{defn}

\begin{rem}
Verdier showed that given a commutative diagram
\begin{equation*}
  \cxymatrix{
    X \ar[r]^u \ar[d]_f & Y \ar[d]^g \ar[r]^{v} & Z \ar[r]^w & \Sigma X \ar[d]^{\Sigma f} \\
    X' \ar[r]_{u'} & Y' \ar[r]_{v'} & Z' \ar[r]_{w'} & \Sigma X' \\
  }
\end{equation*}
in which the rows are exact, there exists a %
fill-in $h : Z \to Z'$ making $(f,g,h)$
into a middling good map of triangles.
We recall the proof and characterize Verdier's extensions in \cref{se:vgood}.
\end{rem}

\begin{defn}\label{def:MappingCone}
The \Def{mapping cone} of a map of triangles~\eqref{eq:map-of-triangles} is the sequence
\[
  \xymatrix@C+2.0pc{
               X' \oplus Y \ar[r]^-{\smat{u' & g \\ 0\; & -v}}
  &            Y' \oplus Z \ar[r]^-{\smat{v' & h \\ 0\; & -w}}
  &            Z' \oplus \Sigma X \ar[r]^-{\smat{w' & \Sigma f \\ 0\; & -\Sigma u}}
  &     \Sigma X' \oplus \Sigma Y ,
  }
\]
where we think of elements of direct sums as column vectors.
The map~\eqref{eq:map-of-triangles} is \Def{good} if the mapping cone is a triangle.
\end{defn}

As an example, the zero map of triangles is always good, since triangles are
closed under direct sum.

In~\cite{Neeman91}*{Theorem~2.3}, Neeman shows that every good map of triangles
is middling good.
He also shows that there are maps which are middling good but not good.
We will prove a slightly stronger claim in \cref{cor:examples}.

\begin{defn}
Let $(f, g, h)$ be a map of candidate triangles as in~\eqref{eq:map-of-triangles} and
let $(f', g', h')$ be another map between the same candidate triangles.
A \Def{chain homotopy} from $(f, g, h)$ to $(f', g', h')$ consists of maps
$F \colon Y \to X'$, $G \colon Z \to Y'$ and $H \colon \Sigma X \to Z'$ such that
$f' - f = Fu + \Sigma^{-1} (w' H)$, $g' - g = Gv + u'F$ and $h' - h = Hw +v'G$, as illustrated in the diagram
\[
\cxymatrix{
X \ar@<-0.5ex>[d] \ar@<0.5ex>[d] \ar[r]^-{u} & Y \ar[dl]_{F} \ar@<-0.5ex>[d] \ar@<0.5ex>[d] \ar[r]^-{v} & Z \ar[dl]_{G} \ar@<-0.5ex>[d] \ar@<0.5ex>[d] \ar[r]^-{w} & \Sigma X \ar[dl]_{H} \ar@<-0.5ex>[d] \ar@<0.5ex>[d] \\
X' \ar[r]_-{u'} & Y' \ar[r]_-{v'} & Z' \ar[r]_-{w'} & \Sigma X'. \\
}
\]
When such a chain homotopy exists, we say that $(f, g, h)$ is \Def{chain homotopic} to $(f', g', h')$.

A map of candidate triangles $(f,g,h)$ is \Def{nullhomotopic} 
if it is chain homotopic to the zero map $(0,0,0)$. 
A candidate triangle $X \to Y \to Z \to \Si X$ is \Def{contractible} if its identity map $(1_X, 1_Y, 1_Z)$ is nullhomotopic (which makes it automatically exact \cite{Neeman01}*{Proposition~1.3.8}). For example, every split triangle is contractible. 
\end{defn}

\begin{rem}\label{rem:good-homotopic}
Neeman observes that if $(f, g, h)$ is chain homotopic to $(f', g', h')$,
then they have isomorphic mapping cones, so are either both good or both not good.
In particular, any map of triangles which is nullhomotopic is good.
It is also immediate that being good is invariant under rotation of triangles.
\end{rem}

\section{Verdier good morphisms}\label{se:vgood}

In this section we introduce the notion of a Verdier good morphism of triangles,
give a characterization of Verdier good morphisms due to Haynes Miller,
and give a simple situation in which good morphisms and Verdier good morphisms agree.

\begin{defn}\label{def:VerdierGood}
A map of triangles as in~\eqref{eq:map-of-triangles} is \Def{Verdier good}
if $h$ can be constructed as in Verdier's proof of the $4 \x 4$ lemma.
This argument is described in the proof of \cite{BeilinsonBD82}*{Proposition~1.1.11}
and summarized in the proof of \cite{Neeman91}*{Lemma~2.1}.
Explicitly, we require that there exist an octahedron for the composite $X \ral{u} Y \ral{g} Y'$: 
\begin{equation*}
\cxymatrix{
X \ar@{=}[d] \ar[r]^-{u} & Y \ar[d]^{g} \ar[r]^-{v} & Z \ar@{-->}[d]^{\al_1} \ar[r]^-{w} & \Si X \ar@{=}[d] \\
X \ar[r]^-{gu} & Y' \ar[d]^{g'} \ar[r]^-{\tild{v}} & A \ar@{-->}[d]^{\be_1} \ar[r]^-{\tild{w}} & \Si X \\
& Y'' \ar[d]^{g''} \ar@{=}[r] & Y'' \ar@{-->}[d]^{\ga_1} & \\
& \Si Y \ar[r]_-{\Si v} & \Si Z, & \\
}
\end{equation*}
and an octahedron for the composite $X \ral{f} X' \ral{u'} Y'$: 
\begin{equation*}
\cxymatrix{
X \ar@{=}[d] \ar[r]^-{f} & X' \ar[d]^{u'} \ar[r]^-{f'} & X'' \ar@{-->}[d]^{\al_2} \ar[r]^-{f''} & \Si X \ar@{=}[d] \\
X \ar[r]^-{u'f = gu} & Y' \ar[d]^{v'} \ar[r]^-{\tild{v}} & A \ar@{-->}[d]^{\be_2} \ar[r]^-{\tild{w}} & \Si X \\
& Z' \ar[d]^{w'} \ar@{=}[r] & Z' \ar@{-->}[d]^{\ga_2} & \\
& \Si X' \ar[r]_-{\Si f'} & \Si X'', & \\
}
\end{equation*}
and that $h \colon Z \to Z'$ is given by $h = \be_2 \circ \al_1$.
\end{defn}

\begin{rem}\label{rem:IsoVerdierGood}
The property of $h$ being a Verdier good fill-in
does not depend on the choice of triangles extending the five maps $u$, $u'$, $f$, $g$, and $gu = u'f$. 
These five triangles may be fixed, and hence were drawn with solid arrows.
Similarly, the notion of Verdier good morphism $(f, g, h)$ is invariant under
pre- and post-composing with isomorphisms of triangles.
One can readily check that the identity morphism of triangles is Verdier good
(by taking $\alpha_1 = 1_Z$ and $\beta_2 = 1_Z$),
so it follows that every isomorphism of triangles is Verdier good.
The notion of good morphism is also invariant under pre- and post-composing
with isomorphisms of triangles, and 
every isomorphism of triangles is good,
since the mapping cone is contractible.  (See~\cite{Neeman91}*{\S 1, Case~1}.)
\end{rem}

\begin{rem}
It is not clear whether being Verdier good is invariant under rotation of triangles.
Since being good is invariant under rotation, our later results that show
that being Verdier good is equivalent to being good in certain situations provide
mild evidence in favour.
It is also not clear whether the notion of Verdier goodness is self-dual.
Given a map of triangles $(f,g,h)$ as in~\eqref{eq:map-of-triangles},
the dual map is $(\Si f, h, g)$, viewed as a map of triangles in $\cTo$
from the bottom row to the top row.
The original map $(f,g,h)$ is Verdier good if and only if the \emph{rotation}
$(\Si g, \Si f, h)$ of the dual map is Verdier good in $\cTo$.
In fact, Verdier goodness is self-dual if and only if Verdier goodness is invariant under rotation of triangles.
\end{rem}

Verdier's argument says that any commutative square admits a Verdier good fill-in. A Verdier good morphism in particular extends to a $4 \x 4$ diagram, but it satisfies more equations. The description in the following lemma was kindly provided to us by Haynes Miller.

\begin{lem}\label{lem:VerdierHaynes}
A map of triangles as in~\eqref{eq:map-of-triangles} is Verdier good if and only if it extends to a $4 \x 4$ diagram
\begin{equation*}%
\cxymatrix{
X \ar[d]_f \ar[r]^-{u} & Y \ar[d]^g \ar[r]^{v} & Z \ar[d]^h \ar[r]^-{w} & \Sigma X \ar[d]^{\Sigma f} \\
X' \ar[d]_{f'} \ar[r]^-{u'} & Y' \ar[d]^{g'} \ar[r]^-{v'} & Z' \ar[d]^{h'} \ar[r]^{w'} & \Sigma X' \ar[d]^{\Si f'} \\
X'' \ar[d]_{f''} \ar@{-->}[r]^-{u''} & Y'' \ar[d]^{g''} \ar@{-->}[r]^-{v''} & Z'' \anti \ar[d]^{h''} \ar@{-->}[r]^-{w''} & \Sigma X'' \ar[d]^{\Sigma f''} \\
\Si X \ar[r]_-{\Si u} & \Si Y \ar[r]_-{\Si v} & \Si Z \ar[r]_-{\Si w} & \Sigma^2 X \\
}
\end{equation*}
(i.e., it is middling good)
and there is an object $A$ together with maps satisfying the following three hexagonal diagrams:
\begin{equation*}%
\vcenter{
\xymatrix @C-1pc @R+0.5pc {
& X \ar[rr]^-{gu=u'f} & & Y' \ar[dl]_{\tild{v}} \ar[dr]^{v'} & \\
Z \circar[ur]^{w} \ar@{-->}[rr]^-{\al_1} & & A \circar[ul]_{\tild{w}} \ar@{-->}[dl]_{\be_1} \ar@{-->}[rr]^-{\be_2} & & Z' \circar[dl]^{(\Si f') w'} \\
& Y'' \circar[ul]^{(\Si v) g''} & & X'', \ar@{-->}[ll]^-{u''} \ar@{-->}[ul]_{\al_2} & \\ 
}
}
\end{equation*}
where each %
triangular shape is either commutative or exact (as indicated by the degree shifts), and each of the three diameters composes to the corresponding vertical map in the $4 \x 4$ diagram; and
\begin{center}
\begin{minipage}{0.25\textwidth}
\begin{equation*}%
\vcenter{
\xymatrix @C-.8pc @R-1pc {
& \Si X \ar[dl]_{\Si f} \ar[dr]^{\Si u} & \\
\Si X' & & \Si Y \\
& A \ar@{-->}[dl]_{\be_2} \ar[uu]_{\tild{w}} \ar@{-->}[dr]^{\be_1} & \\
Z' \ar[uu]^{w'} \ar@{-->}[dr]_{h'} & & Y'' \ar@{-->}[dl]^{v''} \ar[uu]_{g''} \\
& Z'' & \\
}
}
\end{equation*}
\end{minipage}\hspace{0.1\textwidth}
\begin{minipage}{0.25\textwidth}
\begin{equation*}%
\vcenter{
\xymatrix @C-.8pc @R-1pc {
& Y' \ar[dd]_{\tild{v}} & \\
Y \ar[dd]_{v} \ar[ur]^{g} & & X' \ar[ul]_{u'} \ar[dd]^{f'} \\
& A & \\
Z \ar@{-->}[ur]^{\al_1} & \boxed{-1} & X'' \ar@{-->}[ul]_{\al_2} \\
& \Si^{-1} Z'', \ar@{-->}[ul]^{\Si^{-1} h''} \ar@{-->}[ur]_{\Si^{-1} w''} & \\
}
}
\end{equation*}
\end{minipage}
\end{center}
where all squares commute, except for the indicated square which anticommutes.
\end{lem}

\begin{proof}
Let us recall Verdier's argument. 
The two octahedra from \cref{def:VerdierGood} produce maps $X'' \ral{\al_2} A \ral{\be_1} Y''$, to which we apply the octahedral axiom: 
\begin{equation*}
\cxymatrix{
X'' \ar@{=}[d] \ar[r]^-{\al_2} & A \ar[d]^{\be_1} \ar[r]^-{\be_2} & Z' \ar@{-->}[d]^{\al_3} \ar[r]^-{\ga_2} & \Si X'' \ar@{=}[d] \\
X'' \ar[r]^-{\be_1 \al_2} & Y'' \ar[d]^{\ga_1} \ar@{-->}[r]^-{v''} & Z'' \ar@{-->}[d]^{\be_3} \ar@{-->}[r]^-{w''} & \Si X'' \\
& \Si Z \ar[d]^{-\Si \al_1} \ar@{=}[r] & \Si Z \ar@{-->}[d]^{\ga_3} & \\
& \Si A \ar[r]_-{\Si \be_2} & \Si Z'. & \\
}
\end{equation*}
A straightforward %
listing of all equations (and exactness conditions) shows that the three octahedra are equivalent to the three hexagonal diagrams plus two of the equations in the $4 \x 4$ diagram. %
The remaining six equations in the $4 \x 4$ diagram follow from the octahedral equations, by Verdier's proof of the $4 \x 4$ lemma.
\end{proof}

\begin{defn}\label{def:Lightning}
Given two solid triangles
\[
\xymatrix @R-0.8pc {
X \ar[r]^-{u} & Y \ar[r]^-{v} & Z \ar[r]^-{w} & \Sigma X \ar@{-->}[dll]_(0.6){\te} \\
X' \ar[r]_-{u'} & Y' \ar[r]_-{v'} & Z' \ar[r]_-{w'} & \Sigma X' , \\
}
\]
a map $h : Z \to Z'$ is called a \Def{lightning flash} if it can be
expressed as a composite $h = v' \te w$ for some map $\te \colon \Si X \to Y'$. 
\end{defn}

The following is proved in \cite{Neeman91}*{Lemma~2.1} but stated differently.

\begin{lem}\label{lem:AddLightning}
Let $(f,g,h)$ be a Verdier good map of triangles and
let $v' \te w \colon Z \to Z'$ be a lightning flash.
Then $(f, g, h + v' \te w)$ is also Verdier good. \qed
\end{lem}

\begin{prop}\label{pr:ContractibleVerdierGood}
Every map of triangles whose source or target is contractible is 
good and 
Verdier good.
\end{prop}

\begin{proof}
Such a map of triangles is nullhomotopic, hence good, by \cref{rem:good-homotopic}.

Let $(f,g,h)$ be the given map of triangles, and let $\ol{h}$ be a Verdier good
fill-in for the square formed using $f$ and $g$.
By \cite{Neeman91}*{Lemma~2.2} (or its dual in case the target is contractible), %
the difference $h - \ol{h}$ is a lightning flash. 
So by \cref{lem:AddLightning}, %
the fill-in $h$ is also Verdier good.
\end{proof}

\section{Obstruction theory}\label{se:obstructions}

In this section, we give a lifting criterion in terms of good fill-ins.
Given a commutative square
\[
  \xymatrix{
    X \ar[r]^f \ar[d]_u & X' \ar[d]^{u'} \\
    Y \ar[r]_g         & Y' , \\
  }
\]
one can look for obstructions to the
existence of a lift
\[
  \xymatrix{
    X \ar[r]^f \ar[d]_u & X' \ar[d]^{u'} \\
    Y \ar[r]_g \ar@{-->}[ur]^{k} & Y' \\
  }
\]
making both triangles commute.
We'll transpose this picture, to make it fit naturally with the
notation for chain homotopies:
\[
  \xymatrix{
    X \ar[r]^u \ar[d]_f & Y \ar[d]^g \ar@{-->}[dl]_k \\
    X' \ar[r]_{u'} & Y'. \\
  }
\]

A naive way to form an obstruction class is to extend the rows to
triangles and consider the possible fill-in maps
\begin{equation}\label{eq:fill-ins}
  \cxymatrix{
    X \ar[r]^u \ar[d]_f & Y \ar[d]^g \ar[r]^{v} & Z \ar@{-->}[d]^h \ar[r]^w & \Sigma X \ar[d]^{\Sigma f} \\
    X' \ar[r]_{u'} & Y' \ar[r]_{v'} & Z' \ar[r]_{w'} & \Sigma X' \\
  }
\end{equation}
making both squares commute.  %
The collection of such maps is a subset of $[Z, Z']$.
Note that if a lift $k : Y \to X'$ exists, then
$v' g = v' u' k = 0$ and $(\Sigma f) w = (\Sigma k) (\Sigma u) w = 0$,
so the zero map $Z \to Z'$ is one of the fill-ins.
One might hope that the converse holds, i.e., if the fill-in
can be chosen to be zero, then the lift $k$ exists, but the
following example shows that this is not true in general.

\begin{ex}\label{ex:naive}
In the derived category of the integers $\DZ$, let $A[k]$ denote
the chain complex having the abelian group $A$ concentrated in degree $k$.
For an integer $n \geq 2$, let $\ep_n \in \Ext^1_{\Z}(\Z/n,\Z) \cong \Z/n$
denote the canonical generator, sitting in a triangle
\[
\cxymatrix{
\Z[0] \ar[r]^-{n} & \Z[0] \ar[r]^-{q_n} & \Z/n[0] \ar[r]^-{\ep_n} & \Z[1],\\
}
\]
where $q_n \colon \Z \to \Z/n$ denotes the quotient map.
In the morphism of triangles
\[
  \cxymatrix{
    \Z \ar[r]^{q_n} \ar[d]_0 & \Z/n \ar[r]^{\ep_n} \ar[d]^{\ep_n} \ar@{-->}[dl]_k & \Z[1] \ar[d]^0 \ar[r]^{-n} & \Z[1] \ar[d]^0 \\
    \Z/n \ar[r]_{\ep_n}  & \Z[1] \ar[r]_{-n}                 & \Z[1] \ar[r]_{-q_n} & \Z/n[1] , \\
  }
\]
the fill-in map is zero, but no lift $k$ exists.
\end{ex}

In~\cite{ChristensenDI04}, an obstruction theory for such lifts is developed,
working in a model category.
In the case of a stable model category, it does specialize to
giving an obstruction lying in the group $[Z, Z']$,
which is zero if and only if the lift exists.
This is possible because the obstructions are produced using
a strictly commutative starting square.
Fill-in maps produced in this way are better than general fill-ins.
We use this observation to give an obstruction class
in any triangulated category.

\begin{thm}[Lifting criterion]\label{th:obstruction}
Given a solid arrow commutative square
\[
  \xymatrix{
    X \ar[r]^u \ar[d]_f & Y \ar[d]^g \ar@{-->}[dl]_k \\
    X' \ar[r]_{u'} & Y' , \\
  }
\]
choose extensions of the rows to cofibre sequences as in~\eqref{eq:fill-ins}.
Then a lift $k$ exists if and only if the zero map $Z \to Z'$ is a \emph{good} fill-in.
\end{thm}

Define the \Def{obstruction class} to be the subset of $[Z, Z']$
consisting of the good fill-in maps.
Then the theorem says that a lift $k$ exists if and only if the
obstruction class contains zero.
Note that~\cite{Neeman91}*{Theorem~1.8} says that the obstruction class is non-empty.

In \cref{cor:lifting} we extend this result to show that Verdier good
maps can also be used to detect lifting.

\begin{proof}
If the lift $k$ exists, then the triple $(f, g, 0)$ is a map of triangles.
Moreover, it is chain homotopic to the zero map, using $k$ as the only non-zero
component of the chain homotopy.
The zero map of triangles is always good, since the coproduct of
triangles is a triangle.
Therefore, our original map $(f, g, 0)$ is good, and so $0$ is in
the obstruction class.

Conversely, suppose that $(f, g, 0)$ is good.
After rotating, the mapping cone of this map is
\[
  \xymatrix@C+2.8pc{
    \Sigma^{-1} Z' \oplus X\,\ar[r]^-{\smat{-\Sigma^{-1} w' & -f \\ 0\; & u}}
  &          \,X' \oplus Y \ar[r]^-{\smat{-u' & -g \\ 0\; & v}}
  &            Y' \oplus Z \ar[r]^-{\smat{-v' &  0 \\ 0\; & w}}
  &            Z' \oplus \Sigma X ,
  }
\]
which we are assuming is a triangle.  
Since the rightmost map is diagonal, this triangle must be isomorphic
to the sum of the two original triangles (with the first rotated).
That is, there is a matrix
\[
  \bmat{rr} a & b \\ c & d \emat
\]
of maps making the following diagram commute:
\begin{equation}\label{eq1}
  \vcenter{
  \xymatrix@C+3pc@R+1pc{
    \Sigma^{-1} Z' \oplus X\,\ar[r]^-{\smat{-\Sigma^{-1} w' & -f \\ 0\; & u}}
  &          \,X' \oplus Y  \ar[r]^-{\smat{-u' & -g \\ 0\; & v}}
  &            Y' \oplus Z  \ar[r]^-{\smat{-v' &  0 \\ 0\; & w}}
  &            Z' \oplus \Sigma X \\
    \Sigma^{-1} Z' \oplus X\,\ar[r]_-{\smat{-\Sigma^{-1} w' & 0 \\ 0\; & u}}
                            \ar@{=}[u]
  &          \,X' \oplus Y  \ar[r]_-{\smat{-u' &  0 \\ 0\; & v}}
                            \ar[u]_{\smat{a & b \\ c & d}}
  &            Y' \oplus Z  \ar[r]_-{\smat{-v' &  0 \\ 0\; & w}}
                            \ar@{=}[u]
  &            Z' \oplus \Sigma X .
                            \ar@{=}[u]
  }
  }
\end{equation}
The left square %
implies that
\begin{equation}\label{eq2}
  b u = -f  \quad\text{and}\quad   d u = u .
\end{equation}
The middle square %
implies that
\begin{equation}\label{eq4}
  -u' b - g d = 0 .
\end{equation}

Consider the following composite of morphisms of triangles:
include the unprimed triangle into the bottom row of~\eqref{eq1},
map upwards using~\eqref{eq1}, and then project onto the unprimed triangle.
This gives a map of triangles
\[
  \xymatrix{
    X \ar[r]^u & Y \ar[r]^v & Z \ar[r]^w & \Sigma X \\
    X \ar[r]_u \ar[u]^1 & Y \ar[r]_v \ar[u]^d & Z \ar[r]_w \ar[u]^1 & \Sigma X \ar[u]^1 ,
  }
\]
which shows that $d$ is an isomorphism.

Now let $k = -b d^{-1}$.
By~\eqref{eq2},
\[
  k u = -b d^{-1} u = -b u = f ,
\]
and by~\eqref{eq4}
\[
  u' k = -u' b d^{-1} = g .
\]
Thus $k$ is the required lift.
\end{proof}

\begin{cor}
In the above situation, 
if $[Z, Z'] = 0$, then a lift always exists.
\end{cor}

Even this special case seems non-trivial.  The proof
relies on a tricky argument due to Neeman, combined with \cref{th:obstruction}.

\begin{proof}
As mentioned above,~\cite{Neeman91}*{Theorem~1.8} shows that there is always
at least one good fill-in map $Z \to Z'$.  
But if all such maps are zero, it follows that the zero map must be
a good fill-in, so \cref{th:obstruction} applies.
\end{proof}

\begin{ex}
Consider the following lifting problem:
\[
  \xymatrix{
    0 \ar[r] \ar[d] & Y \ar[d]^g \ar@{-->}[dl]_k \\
    X' \ar[r]_{u'} & Y' . \\
  }
\]
In this case, there is a unique fill-in $h = v' g$, which is therefore good.
So we see that $g$ lifts through $u'$ if and only if $v' g = 0$, which
recovers the usual exactness property for mapping into a triangle.
When we take $Y' = 0$ instead of $X = 0$, we get a dual result.
\end{ex}

%
%

%
%
%

\section{Inclusions and projections}\label{se:incproj}

In this section, we give more situations in which good morphisms and
Verdier good morphisms agree.
The situations involve maps of triangles in which some of the component morphisms
are split inclusions or split projections.
We also show that every map of triangles is a composite of
two maps that are good and Verdier good.

\begin{prop}\label{pr:IncUpperTriangular}
\begin{enumerate}
\item \label{item:Inc} For a map of triangles of the form
\[
\xymatrix @C+0.5pc {
X \ar[d]_{\inc} \ar[r]^-{u} & Y \ar[d]^{\inc} \ar[r]^-{v} & Z \ar[d]^{h} \ar[r]^-{w} & \Sigma X \ar[d]^{\inc} \\
X \op X'' \ar[r]_-{\smat{u & c\;\; \\ 0 & \,u''}} & Y \op Y'' \ar[r]_-{\smat{v_1' & v_2'}} & Z' \ar[r]_-{\smat{w_1' \\ w_2'}} & \Sigma X \op \Si X'' \\
}
\]
the following are equivalent:
 \begin{enumerate}
 \item \label{item:IncGood} The map is good.
 \item \label{item:IncVerdier} The map is Verdier good.
 \item \label{item:IncExact} The fill-in $h \colon Z \to Z'$ makes the following diagram exact:
\[
\xymatrix{
X'' \ar[r]^-{\smat{vc\; \\ u''}} & Z \op Y'' \ar[r]^-{\smat{h & v_2'}} & Z' \ar[r]^-{w_2'} & \Si X''. \\
}
\]
 \end{enumerate}
When the matrix entry $c$ is zero (so that the second row is isomorphic to
a direct sum of triangles), these three conditions hold for any $h$.

\item \label{item:Proj} Similarly, a map of triangles of the form
\[
\xymatrix @C+0.5pc {
X \op X'' \ar[d]_{\proj} \ar[r]^-{\smat{u & c\;\; \\ 0 & \,u''}} & Y \op Y'' \ar[d]^{\proj} \ar[r]^-{\smat{v_1' & v_2'}} & Z' \ar[d]^{h} \ar[r]^-{\smat{w_1' \\ w_2'}} & \Sigma X \op \Si X'' \ar[d]^{\proj} \\
X'' \ar[r]_-{u''} & Y'' \ar[r]_-{v''} & Z'' \ar[r]_-{w''} & \Sigma X'' \\
}
\]
is good if and only if it is Verdier good if and only if the fill-in $h \colon Z' \to Z''$ makes the following diagram exact:
\[
\xymatrix @C+1.8pc {
Z' \ar[r]^-{\smat{w_1' \\ h}} & \Si X \op Z'' \ar[r]^-{\smat{\Si u & (\Si c) w''}} & \Si Y \ar[r]^-{\Si v_1'} & \Si Z', \\
}
\]
and these conditions hold when $c$ is zero.
\end{enumerate}
\end{prop}

\begin{proof}
We will prove the statement involving inclusions.
The statement involving projections is proved similarly. %

We begin by showing that conditions~\Eqref{item:IncGood} and \Eqref{item:IncExact} are equivalent.
A straightforward computation with shears shows that the mapping cone
\[
\xymatrix @C+2.1pc @R+0.6pc {
(X \op X'') \op Y \ar[r]^-{\smat{u & c & 1 \\ 0 & u''\!\!\! & 0 \\ 0 & 0 & -v}} & (Y \op Y'') \op Z \ar[r]^-{\smat{v_1' & v_2' & h \\ 0\; & 0\; & -w}} & Z' \op \Si X \ar[r]^-{\smat{w_1' & 1 \\ w_2' & 0 \\ 0\; & - \Si u}} & (\Sigma X \op \Si X'') \op \Si Y \\
}
\]
is isomorphic to the direct sum of the three (candidate) triangles
\begin{equation}\label{eq:3Summands}
\vcenter{
\xymatrix @R-1.1pc {
X \ar[r] & 0 \ar[r] & \Si X \ar[r]^{1} & \Si X \\
Y \ar[r]^{1} & Y \ar[r] & 0 \ar[r] & \Si Y \\
X'' \ar[r]^-{\smat{u'' \\ vc\;}} & Y'' \op Z \ar[r]^-{\smat{v_2' & h}} & Z' \ar[r]^-{w_2'} & \Si X''. \\
}
}
\end{equation}
Therefore, the mapping cone is exact if and only if the bottom row of~\eqref{eq:3Summands} is exact.

We next prove that~\Eqref{item:IncVerdier} implies \Eqref{item:IncExact}.
An octahedron for the composite $X \ral{u} Y \ral{\inc} Y \op Y''$
is of the form:
\[
\xymatrix{
X \ar@{=}[d] \ar[r]^-{u} & Y \ar[d]^{\inc} \ar[r]^-{v} & Z \ar@{-->}[d]^{\al_1 = \smat{\phy\; \\ a_1}} \ar[r]^-{w} & \Si X \ar@{=}[d] \\
X \ar[r]^-{\smat{u \\ 0}} & Y \op Y'' \ar[d]^{\proj} \ar[r]^-{\smat{v & 0 \\ 0 & 1}} & Z \op Y'' \ar@{-->}[d]^{\be_1 = \smat{b_1 & 1}} \ar[r]^-{\smat{w & 0}} & \Si X \\
& Y'' \ar[d]^{0} \ar@{=}[r] & Y'' \ar@{-->}[d]^{\ga_1 = 0} & \\
& \Si Y \ar[r]_-{\Si v} & \Si Z, & \\
}
\]
satisfying the equations
\[
\phy v = v, \quad  a_1 v = 0, \quad w \phy = w \quad\text{and}\quad b_1 v = 0.
\]
The equations $\phy v = v$ and $w \phy = w$ force $\phy$ to be an isomorphism. The third column is exact if and only if $b_1 \phy + a_1 = 0$ holds. Note that picking $\phy = 1_Z$, $a_1 = 0$, and $b_1 = 0$ yields a valid choice.
 
An octahedron for the composite $X \ral{\inc} X \op X'' \ral{u'} Y \op Y''$ is of the form: 
\begin{equation*}
\vcenter{
\xymatrix @R+0.7pc {
X \ar@{=}[d] \ar[r]^-{\inc} & X \op X'' \ar[d]^{\smat{u & c\;\; \\ 0 & u''}} \ar[r]^-{\proj} & X'' \ar@{-->}[d]^{\al_2 = \smat{vc\; \\ u''}} \ar[r]^-{0} & \Si X \ar@{=}[d] \\
X \ar[r]^-{\smat{u \\ 0}} & Y \op Y'' \ar[d]^{\smat{v_1' & v_2'}} \ar[r]^-{\smat{v & 0 \\ 0 & 1}} & Z \op Y'' \ar@{-->}[d]^{\be_2 = \smat{b_2 & v_2'}} \ar[r]^-{\smat{w & 0}} & \Si X \\
& Z' \ar[d]^{\smat{w_1' \\ w_2'}} \ar@{=}[r] & Z' \ar@{-->}[d]^{\ga_2 = w_2'} & \\
& \Si X \op \Si X'' \ar[r]_-{\proj} & \Si X'', & \\
}
}
\end{equation*}
satisfying the equations
\[
b_2 v = v_1', \quad w_1' b_2 = w \quad\text{and}\quad w_2' b_2 = 0.
\]
Those are precisely the equations for $b_2 \colon Z \to Z'$ being a fill-in. Verdier good fill-ins are those of the form
\[
\be_2 \circ \al_1 = b_2 \phy + v_2' a.
\]
Using the equations $\phy v = v$ and $a_1 v = 0$, we obtain isomorphisms of rows
\[
\xymatrix @R+0.4pc @C+1.8pc {
X'' \ar@{=}[d] \ar[r]^-{\smat{vc\; \\ u''}} & Z \op Y'' \ar[d]_{\cong}^{\smat{\phy^{-1} & 0 \\ 0\;\;\;\; & 1}} \ar[r]^-{\smat{b_2 & v_2'}} & Z' \ar@{=}[d] \ar[r]^-{w_2'} & \Si X'' \ar@{=}[d] \\
X'' \ar@{=}[d] \ar[r]^-{\smat{vc\; \\ u''}} & Z \op Y'' \ar[d]_{\cong}^{\smat{1 & 0 \\ -a & 1}} \ar[r]^-{\smat{b_2 \phy & v_2'}} & Z' \ar@{=}[d] \ar[r]^-{w_2'} & \Si X'' \ar@{=}[d] \\
X'' \ar[r]_-{\smat{vc\; \\ u''}} & Z \op Y'' \ar[r]_-{\smat{b_2 \phy + v_2' a & v_2'}} & Z' \ar[r]_-{w_2'} & \Si X''. \\
}
\]
Since the top row is exact, the bottom row is exact, proving the exactness property of Verdier good fill-ins.

To see that~\Eqref{item:IncExact} implies \Eqref{item:IncVerdier}, let $h \colon Z \to Z'$ be a fill-in satisfying~\Eqref{item:IncExact}.
Then the choices $\phy = 1_Z$, $a=0$, $b_1 = 0$, and $b_2 = h$ exhibit $h$ as Verdier good.

When $c = 0$, we verify condition~\Eqref{item:IncExact}.
Without loss of generality, we may assume that the bottom row is a direct sum
of triangles, so the map we are considering is of the form
\[
\xymatrix @C+0.5pc {
X \ar[d]_{\inc} \ar[r]^-{u} & Y \ar[d]^{\inc} \ar[r]^-{v} & Z \ar[d]^{h = \smat{h_1 \\ h_2}} \ar[r]^-{w} & \Sigma X \ar[d]^{\inc} \\
X \op X'' \ar[r]_-{\smat{u & 0\;\; \\ 0 & u''}} & Y \op Y'' \ar[r]_-{\smat{v & 0\;\; \\ 0 & v''}} & Z \op Z'' \ar[r]_-{\smat{w & 0\;\; \\ 0 & w''}} & \Sigma X \op \Si X'' . \\
}
\]
Note that $h_1$ is necessarily an isomorphism and $w'' h_2 = 0$.
The following isomorphism of triangles
\[
\xymatrix @C+0.9pc @R+0.4pc {
X'' \ar@{=}[d] \ar[r]^-{\smat{u'' \\ 0\;\;}} & Y'' \op Z \ar@{=}[d] \ar[r]^-{\smat{0\;\; & h_1 \\ v'' & h_2}} & Z \op Z'' \ar[r]^-{\smat{0 & w''}} & \Sigma X'' \ar@{=}[d] \\
X'' \ar[r]_-{\smat{u'' \\ 0\;\;}} & Y'' \op Z \ar[r]_-{\smat{0\;\; & 1 \\ v'' & 0}} & Z \op Z'' \ar[r]_-{\smat{0 & w''}} \ar[u]^{\cong}_-{\smat{h_1 & 0 \\ h_2 & 1}} & \Sigma X'' \\
}
\]
shows that the top row is isomorphic to a direct sum of exact triangles, and is
therefore exact, proving~\Eqref{item:IncExact}.
\end{proof}

\begin{lem}\label{lem:IncMappingCone}
Let $(f,g,h)$ be a good map of triangles.
\begin{enumerate}
 \item \label{item:IncMappingCone} The inclusion of the target of $(f,g,h)$ into its mapping cone is a good map of triangles.
 \item \label{item:ProjMappingCone} The projection
\[
\xymatrix @C+0.7pc {
X' \op Y \ar[d]_{\proj} \ar[r]^-{\smat{u' & g \\ 0\; & -v}} & Y' \op Z \ar[d]^{\proj} \ar[r]^-{\smat{v' & h \\ 0\; & -w}} & Z' \op \Si X \ar[d]^{\proj} \ar[r]^-{\smat{w' & \Si f \\ 0\; & -\Si u}} & \Si X' \op \Si Y \ar[d]^{\proj} \\
Y \ar[r]_-{-v} & Z \ar[r]_-{-w} & \Si X \ar[r]_-{- \Si u} & \Si Y. \\
}
\] 
from the mapping cone onto the translation of the source of $(f,g,h)$ is a good map of triangles.
\end{enumerate}
\end{lem}

\begin{proof}
We will prove the first statement; the second statement is dual. 
Using the notation from~\eqref{eq:map-of-triangles}, the inclusion of the target of $(f,g,h)$ into its mapping cone is
\[
\xymatrix @C+0.5pc {
X' \ar[d]_{\inc} \ar[r]^-{u'} & Y' \ar[d]^{\inc} \ar[r]^-{v'} & Z' \ar[d]^{\inc} \ar[r]^-{w'} & \Sigma X' \ar[d]^{\inc} \\
X' \op Y \ar[r]_-{\smat{u' & g \\ 0\; & -v}} & Y' \op Z \ar[r]_-{\smat{v' & h \\ 0\; & -w}} & Z' \op \Si X \ar[r]_-{\smat{w' & \Si f \\ 0\; & -\Si u}} & \Sigma X' \op \Si Y. \\
}
\]
By \cref{pr:IncUpperTriangular}, that map is good if and only if the candidate triangle
\[
\xymatrix @C+0.5pc {
Y \ar[r]^-{\smat{v' g \\ -v}} & Z' \op Z \ar[r]^-{\smat{1 & h \\ 0 & -w}} & Z' \op \Si X \ar[r]^-{\smat{0 & - \Si u}} & \Si Y \\
}
\]
is exact. By the isomorphism of rows
\[
\xymatrix @C+0.5pc {
Y \ar@{=}[d] \ar[r]^-{\smat{v' g \\ -v}} & Z' \op Z \ar[d]_{\cong}^{\smat{1 & h \\ 0 & 1}} \ar[r]^-{\smat{1 & h \\ 0 & -w}} & Z' \op \Si X \ar@{=}[d] \ar[r]^-{\smat{0 & - \Si u}} & \Si Y \ar@{=}[d] \\
Y \ar[r]_-{\smat{0 \\ -v}} & Z' \op Z \ar[r]_-{\smat{1 & 0 \\ 0 & -w}} & Z' \op \Si X \ar[r]_-{\smat{0 & - \Si u}} & \Si Y, \\
}
\]
the top row is indeed exact.
\end{proof}

\begin{lem}\label{lem:IntoProduct}
Let $(f,g,h)$ be a map of triangles as in~\eqref{eq:map-of-triangles} and $(\ol{f}, \ol{g}, \ol{h})$ a map of triangles from the same source triangle, as in the diagram
\[
\cxymatrix{
X \ar[d]_{\ol{f}} \ar[r]^-{u} & Y \ar[d]^{\ol{g}} \ar[r]^-{v} & Z \ar[d]^{\ol{h}} \ar[r]^-{w} & \Sigma X \ar[d]^{\Si \ol{f}} \\
\ol{X} \ar[r]_-{\ol{u}} & \ol{Y} \ar[r]_-{\ol{v}} & \ol{Z} \ar[r]_-{\ol{w}} & \Sigma \ol{X}. \\
}
\]
Assume that one of the two maps is nullhomotopic. Then both maps are good if and only if the map of triangles
\[
\cxymatrix{
X \ar[d]_{\smat{f \\ \ol{f}}} \ar[r]^-{u} & Y \ar[d]^{\smat{g \\ \ol{g}}} \ar[r]^-{v} & Z \ar[d]^{\smat{h \\ \ol{h}}} \ar[r]^-{w} & \Sigma X \ar[d]^{\smat{\Si f \\ \Si \ol{f}}} \\
X' \op \ol{X} \ar[r]_-{\smat{u' & 0 \\ 0\; & \ol{u}}} & Y' \op \ol{Y} \ar[r]_-{\smat{v' & 0 \\ 0\; & \ol{v}}} & Z' \op \ol{Z} \ar[r]_-{\smat{w' & 0 \\ 0\; & \ol{w}}} & \Si X' \op \Sigma \ol{X} \\
}
\]
is good. 

The analogous statement for two maps of triangles to the same target triangle holds.
\end{lem}

\begin{proof}
Assume that the map $(\ol{f}, \ol{g}, \ol{h})$ is nullhomotopic (in particular good). Then the map \linebreak $(\smat{f \\ \ol{f}}, \smat{g \\ \ol{g}}, \smat{h \\ \ol{h}})$ is chain homotopic to $(\smat{f \\ 0}, \smat{g \\ 0}, \smat{h \\ 0})$. The latter has as mapping cone the mapping cone of $(f,g,h)$ direct sum with $\ol{X} \ral{\ol{u}} \ol{Y} \ral{\ol{v}} \ol{Z} \ral{\ol{w}} \Si \ol{X}$. %
\end{proof}

The following uses a variation on an argument of Neeman \cite{Neeman91}*{before Remark~1.10}.

\begin{prop}\label{pr:compositeVgood}
Every map of triangles is a composite of two maps that are good and Verdier good.
\end{prop}

In particular, a composite of Verdier good maps need not be Verdier good, and a composite of middling good maps need not be middling good.

\begin{proof}
Let $(f,g,h)$ be a map of triangles as in~\eqref{eq:map-of-triangles}. Consider the commutative diagram
\[
\xymatrix @R+0.2pc @C+0.2pc {
X \ar[d]_{\inc} \ar[r]^-{u} & Y \ar[d]^{\inc} \ar[r]^-{v} & Z \ar[d]^{\inc} \ar[r]^-{w} & \Sigma X \ar[d]^{\inc} \\
X \op X' \ar[d]_{\smat{1 & 0 \\ f & 1}}^{\cong} \ar[r]^-{\smat{u & 0\; \\ 0 & u'}} & Y \op Y' \ar[d]^{\smat{1 & 0 \\ g & 1}}_{\cong} \ar[r]^-{\smat{v & 0\; \\ 0 & v'}} & Z \op Z' \ar[d]^{\smat{1 & 0 \\ h & 1}}_{\cong} \ar[r]^-{\smat{w & 0\; \\ 0 & w'}} & \Sigma X \op \Si X' \ar[d]^{\smat{1 & 0 \\ \Si f & 1}}_{\cong} \\
X \op X' \ar[d]_{\proj} \ar[r]_-{\smat{u & 0\; \\ 0 & u'}} & Y \op Y' \ar[d]^{\proj} \ar[r]_-{\smat{v & 0\; \\ 0 & v'}} & Z \op Z' \ar[d]^{\proj} \ar[r]_-{\smat{w & 0\; \\ 0 & w'}} & \Sigma X \op \Si X' \ar[d]^{\proj} \\
X' \ar[r]_-{u'} & Y' \ar[r]_-{v'} & Z' \ar[r]_-{w'} & \Sigma X', \\
}
\]
which consists of three maps of triangles whose composite is $(f,g,h)$. By \cref{pr:IncUpperTriangular}, the top and the bottom maps are good and Verdier good (as is the middle map, by \cref{rem:IsoVerdierGood}). By \cref{rem:IsoVerdierGood}, the composite of the top two (or bottom two) maps is good and Verdier good.
\end{proof}

\section{Homotopy cartesian squares}\label{se:hocart}

In this section, we recall the notion of a homotopy cartesian square and
build on work of Neeman showing a relationship between such squares and
good morphisms of the form $(1, g, h)$.
We then prove a pasting lemma for homotopy cartesian squares.
We also characterize the Verdier good morphisms of the form $(1,g,h)$,
and deduce that every good morphism of the form $(1,g,h)$ is Verdier good.

We will use \cite{Neeman01}*{Definition~1.4.1, Lemma~1.4.3}, but different preprint versions of the book have different sign conventions. Let us recall the sign convention and check the details of the proof.

\begin{defn}\label{def:HoCart}
A commutative square
\begin{equation}\label{eq:HoCart}
\cxymatrix{
Y \ar[d]_g \ar[r]^-{f} & Z \ar[d]^{g'} \\
Y' \ar[r]_-{f'} & Z' \\
}
\end{equation}
is called \Def{homotopy cartesian} if there is an exact triangle 
\begin{equation}\label{eq:MayerVietoris}
\xymatrix @C+1.5pc {
Y \ar[r]^-{\smat{g \\ f}} & Y' \op Z \ar[r]^-{\smat{f' & -g'}} & Z' \ar[r]^-{\del} & \Si Y \\
}
\end{equation}
for some map $\del \colon Z' \to \Si Y$, called the \Def{differential}. 
\end{defn}
Note that exactness of \eqref{eq:MayerVietoris} is equivalent to exactness of
\[
\xymatrix @C+1.5pc {
Y \ar[r]^-{\smat{g \\ -f}} & Y' \op Z \ar[r]^-{\smat{f' & g'}} & Z' \ar[r]^-{\del} & \Si Y, \\
}
\]
via the automorphism $-1_Z \colon Z \to Z$.

\begin{rem}\label{rem:HoCartFlip}
The horizontal and vertical directions play asymmetric roles in the sign convention. The square \eqref{eq:HoCart} is homotopy cartesian with differential $\del \colon Z' \to \Si Y$ if and only if its diagonal reflection
\[
\cxymatrix{
Y \ar[d]_f \ar[r]^-{g} & Y' \ar[d]^{f'} \\
Z \ar[r]_-{g'} & Z' \\
}
\]
is homotopy cartesian with differential $-\del \colon Z' \to \Si Y$.
\end{rem}

The following is stated in \cite{Neeman01}*{Remark~1.4.5}. Here we fill in some details of the proof.

\begin{lem}\label{lem:HoCart}
Consider a homotopy cartesian square
\[
\cxymatrix{
Y \ar[d]_g \ar[r]^-{v} & Z \ar[d]^h \\
Y' \ar[r]_-{v'} & Z' \\
}
\]
with given differential $\del \colon Z' \to \Si Y$. 
Then the square extends to %
a good map of triangles of the form $(1,g,h)$, as illustrated in the diagram
\[
\cxymatrix{
X \ar@{=}[d] \ar[r]^-{u} & Y \ar[d]_g \ar[r]^-{v} & Z \ar[d]^h \ar[r]^-{w} & \Si X \ar@{=}[d] \\
X \ar[r]_-{u'} & Y' \ar[r]_-{v'} & Z' \ar[r]_-{w'} & \Si X, \\
}
\]
satisfying $\del = (\Si u) w'$.

Moreover, we may prescribe the top row or the bottom row.
\end{lem}

\begin{proof}
The map is constructed in \cite{Neeman01}*{Lemma~1.4.4}. Let us recall Neeman's proof and check that the resulting map is good. 

Assume given a triangle $X \ral{u} Y \ral{v} Z \ral{w} \Si X$, prescribed as top row. Let $\phy \colon Z' \to \Si X$ be a good fill-in in the diagram
\begin{equation}\label{eq:GoodFillin}
\cxymatrix{
Y \ar@{=}[d] \ar[r]^-{\smat{g \\ v}} & Y' \op Z \ar[d]^{\proj} \ar[r]^-{\smat{v' & -h}} & Z' \ar@{-->}[d]^{\phy} \ar[r]^{\del} & \Si Y \ar@{=}[d] \\
Y \ar[r]_-{v} & Z \ar[r]_-{w} & \Si X \ar[r]_-{-\Si u} & \Si Y. \\
}
\end{equation}
A straightforward computation shows that its mapping cone is isomorphic to the direct sum of the three summands
\[
\xymatrix @R-1.1pc @C+1.1pc {
Z \ar[r]^-{1} & Z \ar[r] & 0 \ar[r] & \Si Z \\
Y \ar[r] & 0 \ar[r] & \Si Y  \ar[r]^-{1} & \Si Y \\
Y' \ar[r]^-{v'} & Z' \ar[r]^-{-\phy} & \Si X \ar[r]^{- (\Si g)(\Si u)} & \Si Y'. \\
}
\]
Via the isomorphism, the inclusion of the bottom row of~\eqref{eq:GoodFillin} into the mapping cone has as component into the third summand
\begin{equation}\label{eq:ThirdComponent}
\vcenter{
\xymatrix @C+1pc {
Y \ar[d]_g \ar[r]^-{v} & Z \ar[d]^h \ar[r]^-{w} & \Si X \ar@{=}[d] \ar[r]^{-\Si u} & \Si Y \ar[d]^{\Si g} \\
Y' \ar[r]_-{v'} & Z' \ar[r]_-{-\phy} & \Si X \ar[r]_{-(\Si g)(\Si u)} & \Si Y', \\
}
}
\end{equation}
which is the desired map (up to rotation), taking $w' := -\phy$ and $u' := gu$. By \cref{lem:IncMappingCone}, the inclusion of the bottom row of~\eqref{eq:GoodFillin} into the mapping cone is good. By \cref{lem:IntoProduct}, its third component \eqref{eq:ThirdComponent} is also good. Furthermore, the differential satisfies $\del = (\Si u) (-\phy) = (\Si u) w'$. 

Now assume given a triangle $X \ral{u'} Y' \ral{v'} Z' \ral{w'} \Si X$, prescribed as bottom row. Applying an appropriate automorphism of $X$ to the map $(1,g,h)$ found above will adjust the bottom row to the prescribed one, modify the top row, and keep the differential $\del$ as it is. 
\end{proof}

We next give a characterization of the good maps of the form $(1,g,h)$.

\begin{defn}
In a candidate triangle $X \ral{u} Y \ral{v} Z \ral{w} \Si X$, the map $u$ is \Def{replaceable} with \Def{replacing map} $\hatt{u}$ if
\[
\cxymatrix{
X \ar[r]^-{\hatt{u}} & Y \ar[r]^-{v} & Z \ar[r]^-{w} & \Si X \\
}
\]
is exact. Replaceability of $v$ or $w$ is defined similarly. The candidate triangle is %
\Def{replaceably exact} 
if its three maps are replaceable \cite{Vaknin01}*{Definition~1.3}.
\end{defn}

\begin{prop}\label{prop:HoCartFillin}
Consider a map of triangles $(1,g,h)$ of the form
\[
\xymatrix{
X \ar@{=}[d] \ar[r]^-{u} & Y \ar[d]^{g} \ar[r]^-{v} & Z \ar[d]^{h} \ar[r]^-{w} & \Si X \ar@{=}[d] \\
X \ar[r]_-{u'} & Y' \ar[r]_-{v'} & Z' \ar[r]_-{w'} & \Si X \\
}
\]
and the candidate triangle
\begin{equation}\label{eq:MV}
\vcenter{
\xymatrix @C+0.5pc{
Y \ar[r]^-{\smat{g \\ v}} & Y' \op Z \ar[r]^-{\smat{v' & -h}} & Z' \ar[r]^-{(\Si u)w'} & \Si Y .
}
}
\end{equation}
\begin{enumerate}
\item \label{item:1ghGood} The map %
is good if and only if the middle square is homotopy cartesian with differential $\del = (\Si u) w' \colon Z' \to \Si Y$ (i.e.,~\eqref{eq:MV} is exact).
\item \label{item:MiddleHoCart} The middle square is homotopy cartesian if and only if the %
candidate triangle~\eqref{eq:MV} is replaceably exact.
\end{enumerate}
\end{prop} 

Note that if the cofibre fill-in $h \colon Z \to Z'$ were missing, we could choose a good one, and if the fibre fill-in $g \colon Y \to Y'$ were missing, we could choose a good one.

\begin{proof} 
Part~\eqref{item:1ghGood} is a slightly stronger formulation of \cite{Neeman01}*{Lemma~1.4.3}, based on the same proof, which we recall here. The mapping cone of $(1,g,h)$ is
\[
\xymatrix @C+1.5pc {
X \op Y \ar[r]^-{\smat{gu & g \\ 0\; & -v}} & Y' \op Z \ar[r]^-{\smat{v' & h \\ 0\; & -w}} & Z' \op \Si X \ar[r]^-{\smat{w' & 1 \\ 0\; & -\Si u}} & \Si X \op \Si Y. \\
}
\]
A straightforward computation shows that it is isomorphic to the direct sum of
the candidate triangle~\eqref{eq:MV} and the triangle
\[
\cxymatrix {
X \ar[r] & 0 \ar[r] & \Si X \ar[r]^1 & \Si X. \\
}
\]

\eqref{item:MiddleHoCart} By definition, the middle square is homotopy cartesian if and only if~\eqref{eq:MV} is replaceable on the right. Assuming that the middle square is homotopy cartesian, we will now show that~\eqref{eq:MV} is replaceable in the middle. Applying \cref{lem:HoCart} to the middle square and prescribing the top row, there is a good map of triangles
\[
\cxymatrix{
X \ar@{=}[d] \ar[r]^-{u} & Y \ar[d]^{g} \ar[r]^-{v} & Z \ar[d]^{h} \ar[r]^-{w} & \Sigma X \ar@{=}[d] \\
X \ar[r]_-{u'} & Y' \ar[r]_-{v'} & Z' \ar[r]_-{\hatt{w}'} & \Sigma X. \\
}
\]
Now postcompose with the isomorphism of triangles
\[
\cxymatrix{
X \ar@{=}[d] \ar[r]^-{u'} & Y' \ar@{=}[d] \ar[r]^-{v'} & Z' \ar@{-->}[d]^{\te}_{\cong} \ar[r]^-{\hatt{w}'} & \Sigma X \ar@{=}[d] \\
X \ar[r]_-{u'} & Y' \ar[r]_-{v'} & Z' \ar[r]_-{w'} & \Sigma X, \\
}
\]
where $\te \colon Z' \to Z'$ is any fill-in. By part~\eqref{item:1ghGood}, %
\[
\xymatrix @C+0.7pc {
Y \ar[r]^-{\smat{g \\ v}} & Y' \op Z \ar[r]^-{\smat{v' & - \te h}} & Z' \ar[r]^-{(\Si u) w'} & \Si Y \\
}
\]
is a %
triangle. 
Dualizing the argument, prescribing the bottom row, there is an automorphism $\Ga \colon Y \ral{\cong} Y$ such that 
\[
\xymatrix @C+0.5pc {
Y \ar[r]^-{\smat{g \Ga \\ v}} & Y' \op Z \ar[r]^-{\smat{v' & -h}} & Z' \ar[r]^-{(\Si u) w'} & \Si Y \\
}
\]
is a %
triangle. This shows that~\eqref{eq:MV} is replaceable on the left.
\end{proof}

We illustrate the usefulness of our results by proving the following fact
about homotopy cartesian squares,
whose statement does not involve the notion of good morphism.

\begin{prop}[Pasting Lemma]\label{pr:Pasting}
Consider a diagram
\begin{equation*}
\cxymatrix{
X \ar[d]_f \ar[r]^-{\phy} & Y \ar[d]^g \ar[r]^-{\psi} & Z \ar[d]^h \\
X' \ar[r]_-{\phy'}  & Y' \ar[r]_-{\psi'} & Z'. \\
}
\end{equation*}
If the two squares are homotopy cartesian, then so is their pasting, the big rectangle.

Moreover, given differentials $\del_L \colon Y' \to \Si X$ and $\del_R \colon Z' \to \Si Y$ of the left square and right square respectively, there exists a differential $\del_P \colon Z' \to \Si X$ for the pasted %
rectangle satisfying
\[
\begin{cases}
(\Si \phy) \del_P = \del_R \\
\del_P \psi' = \del_L. \\
\end{cases}
\]
\end{prop}

\begin{proof}
We claim that the square
\[
\cxymatrix{
X \ar[d]_{\smat{f\; \\ \psi \phy}} \ar[r]^-{\phy} & Y \ar[d]^{\smat{g \\ \psi}} \\ 
X' \op Z \ar[r]_-{\smat{\phy' & 0 \\ 0\; & 1}} & Y' \op Z \\ 
}
\]
is homotopy cartesian with differential $\smat{\del_L & 0}$.
This follows from the commutative diagram
\[
\vcenter{
\xymatrix @C+0.9pc {
X \ar@{=}[d] \ar[r]^-{\smat{f \\ 0 \\ \phy}} & X' \op Z \op Y \ar[d]^s_{\cong} \ar[r]^-{\smat{\psi' & 0 & -g \\ 0\; & 1 & 0 }} & Y' \op Z \ar@{=}[d] \ar[r]^--{\smat{\del_L & 0}} & \Si X \ar@{=}[d] \\ 
X \ar[r]_-{\smat{f\; \\ \psi \phy \\ \phy\;}} & X' \op Z \op Y \ar[r]_-{\smat{\psi' & 0 & -g \\ 0\; & 1 & -\psi }} & Y' \op Z \ar[r]_-{\smat{\del_L & 0}} & \Si X \\ 
}
}
\]
involving the shear isomorphism $s$ with matrix
\[
  \bmat{rrr} 1 & 0 & 0 \\ 0 & 1 & \psi \\ 0 & 0 & 1 \emat \raisebox{-13pt}{.}
\]
The top row is exact since it is a direct sum of the triangle that
shows that the left square is homotopy cartesian with a trivial triangle
involving the identity map on $Z$.
Therefore, the bottom row is exact, showing that the square is homotopy cartesian.

Keeping that specific choice of differential and applying \cref{lem:HoCart}, there exists a (good) map of triangles 
\begin{equation}\label{eq:HocartCombined}
\vcenter{
\xymatrix @C+1.2pc {
X \ar[d]_{\phy} \ar[r]^-{\smat{f\; \\ \psi \phy}} & X' \op Z \ar[d]^{\smat{\phy' & 0 \\ 0\; & 1}} \ar[r]^-{\smat{\psi' \phy' & -h}} & Z' \ar@{=}[d] \ar[r]^-{\del_P} & \Si X \ar[d]^{\Si \phy} \\ 
Y \ar[r]_-{\smat{g \\ \psi}} & Y' \op Z \ar[r]_-{\smat{\psi' & -h}} & Z' \ar[r]_-{\del_R} & \Si Y \\ 
}
}
\end{equation}
satisfying
\[
\smat{-\del_L & 0} = -\del_P \smat{\psi' & -h} = \smat{-\del_P \psi' & \; 0}.
\]
(The minus sign on the left-hand side comes from the reflection, as in \cref{rem:HoCartFlip}.) 
Here we prescribed the bottom row, coming from the fact that the original right square was homotopy cartesian with given differential $\del_R \colon Z' \to \Si Y$. The top row of~\eqref{eq:HocartCombined} exhibits the big rectangle as homotopy cartesian with differential $\del_P \colon Z' \to \Si X$. The chosen differential $\del_P$ satisfies $\del_R = (\Si \phy) \del_P$ and $\del_L = \del_P \psi'$, as desired.
\end{proof}

\begin{rem}
\cite{LinZ19}*{Lemma~2.7} is a form of pasting lemma. However, in the case $n=3$ (which is the case we are considering, of usual triangles), their statement only allows pasting an isomorphism of arrows, i.e., our diagram with $\psi$ and $\psi'$ being isomorphisms.  
\end{rem}

\begin{cor}
Consider a composite of good maps of triangles of the form
\[
\cxymatrix{
X \ar@{=}[d] \ar[r]^-{u_1} & Y_1 \ar[d]^{g_1} \ar[r]^-{v_1} & Z_1 \ar[d]^{h_1} \ar[r]^-{w_1} & \Sigma X \ar@{=}[d] \\
X \ar@{=}[d] \ar[r]^-{u_2} & Y_2 \ar[d]^{g_2} \ar[r]^-{v_2} & Z_1 \ar[d]^{h_2} \ar[r]^-{w_2} & \Sigma X \ar@{=}[d] \\
X \ar[r]_-{u_3} & Y_3 \ar[r]_-{v_3} & Z_3 \ar[r]_-{w_3} & \Sigma X. \\
}
\]
Then the %
candidate triangle
\begin{equation}\label{eq:MVPasting}
\xymatrix @C+1.2pc {
Y_1 \ar[r]^-{\smat{g_2 g_1 \\ v_1}} & Y_3 \op Z_1 \ar[r]^-{\smat{v_3 & -h_2 h_1}} & Z_3 \ar[r]^-{(\Si u_1) w_3} & \Si Y_1 \\
}
\end{equation}
is replaceably exact.
\end{cor}

\begin{proof}
By \cref{pr:Pasting}, the pasting of the top middle and bottom middle squares is homotopy cartesian. The claim then follows from \cref{prop:HoCartFillin}~\eqref{item:MiddleHoCart}.
\end{proof}

One might hope that the differential for the pasted rectangle could be chosen
to be $(\Si u_1) w_3$, or, in other words, that the candidate triangle~\eqref{eq:MVPasting}
is exact. 
In \cref{ex:F2C4composite}, we will show that this need \emph{not} be the case. 
In \cite{Vaknin01}*{after Proposition~1.16}, Vaknin shows that~\eqref{eq:MVPasting} is ``virtual'',
which is weaker than being replaceably exact.
\medskip

We now give a characterization of when a map $(1,g,h)$
of triangles is Verdier good.

\begin{prop}\label{pr:1ghVerdier}
A map of triangles of the form %
\[
\cxymatrix{
X \ar@{=}[d] \ar[r]^-{u} & Y \ar[d]^{g} \ar[r]^-{v} & Z \ar[d]^{h} \ar[r]^-{w} & \Sigma X \ar@{=}[d] \\
X \ar[r]_-{u'} & Y' \ar[r]_-{v'} & Z' \ar[r]_-{w'} & \Sigma X \\
}
\]
is Verdier good if and only if it extends to an octahedron, i.e., it appears as the top two rows of an octahedron.
\end{prop}

\begin{proof}
An octahedron for the composite $X \ral{u} Y \ral{g} Y'$ is of the form:
\begin{equation}\label{eq:OctahedronUG1}
\cxymatrix{
X \ar@{=}[d] \ar[r]^-{u} & Y \ar[d]^{g} \ar[r]^-{v} & Z \ar@{-->}[d]^{\al_1} \ar[r]^-{w} & \Si X \ar@{=}[d] \\
X \ar[r]^-{gu = u'} & Y' \ar[d]^{g'} \ar[r]^-{v'} & Z' \ar@{-->}[d]^{\be_1} \ar[r]^-{w'} & \Si X \\
& Y'' \ar[d]^{g''} \ar@{=}[r] & Y'' \ar@{-->}[d]^{\ga_1} & \\
& \Si Y \ar[r]_-{\Si v} & \Si Z. & \\
}
\end{equation}
An octahedron for the composite $X \ral{1} X \ral{u'} Y'$ is of the form:
\[
\cxymatrix{
X \ar@{=}[d] \ar[r]^-{1} & X \ar[d]^{u'} \ar[r] & 0 \ar[d] \ar[r] & \Si X \ar@{=}[d] \\
X \ar[r]^-{u' = gu} & Y' \ar[d]^{v'} \ar[r]^-{v'} & Z' \ar@{-->}[d]^{\be_2}_{\cong} \ar[r]^-{w'} & \Si X \\
& Z' \ar[d]^{w'} \ar@{=}[r] & Z' \ar[d] & \\
& \Si X \ar[r] & 0, & \\
}
\]
with the automorphism $\be_2$ satisfying $\be_2 v' = v'$ %
and $w' \be_2 = w'$. %

($\impliedby$) If the map of triangles $(1,g,h)$ extends to an octahedron, then that octahedron is a valid choice for the octahedron~\eqref{eq:OctahedronUG1}, in particular with $\al_1 = h$. Picking $\be_2 = 1_{Z'}$ yields a valid choice for the second octahedron, exhibiting $h = 1_{Z'} h$ as Verdier good.

($\implies$) Assume that $h$ is Verdier good, exhibited by two octahedra as above. Using the automorphism $\be_2 \colon Z' \ral{\cong} Z'$ to modify the first octahedron~\eqref{eq:OctahedronUG1} yields an octahedron
\[
\xymatrix @C+0.9pc {
X \ar@{=}[d] \ar[r]^-{u} & Y \ar[d]^{g} \ar[r]^-{v} & Z \ar[d]^{\be_2 \al_1 = h} \ar[r]^-{w} & \Si X \ar@{=}[d] \\
X \ar[r]^-{u'} & Y' \ar[d]^{g'} \ar[r]^-{\be_2 v' = v'} & Z' \ar[d]^{\be_1 \be_2^{-1}} \ar[r]^-{w' \be_2^{-1} = w'} & \Si X \\
& Y'' \ar[d]^{g''} \ar@{=}[r] & Y'' \ar@{-->}[d]^{\ga_1} & \\
& \Si Y \ar[r]_-{\Si v} & \Si Z & \\
}
\]
that extends the map $(1,g,h)$.
\end{proof}

The proof of \cite{Neeman01}*{Proposition~1.4.6} then implies the following.

\begin{cor}\label{cor:1ghVerdier}
If a map of triangles of the form $(1,g,h)$ is good, then it is Verdier good.  \qed
\end{cor}

We don't know whether Verdier good implies good for such maps of triangles.
The discussion in \cite{May01}*{Remark~3.7} suggests that this is \emph{not} the case, without providing an explicit counterexample.

The following examples illustrate the difference between a map $(1,g,h)$ being good and merely having a homotopy cartesian middle square. They also show that we may not always prescribe both rows in \cref{lem:HoCart}. 
Furthermore, they provide examples of maps of the form $(1,g,h)$ that are middling good but not Verdier good.

\begin{ex}\label{ex:F2C4}
Let $C_4$ be the cyclic group of order $4$, with generator $g$, and consider the group algebra $R = \F_2 C_4 \cong \F_2[x]/x^4$, with $x := g-1$. We will work in the stable module category $\StMod(R)$, as in \cite{ChristensenFrankland17}*{Appendix~A}. See also \cite{Benson98rep1}*{\S 2.1} %
or \cite{Carlson96}*{\S 5} %
for background on stable module categories. %
For $r \in R$, let $\mu_r \colon R/x^i \to R/x^j$ denote the $R$-module map sending $1$ to $r$, when it is defined. 
In the map of triangles
\begin{equation}\label{eq:F2C4good}
\cxymatrix{
R/x^3 \op R/x^2 \ar@{=}[d] \ar[r]^-{\smat{1 & 0 \\ 0 & \mu_x}} & R/x^3 \op R/x^3 \ar[d]^{\smat{\mu_1 & 0}} \ar[r]^-{\smat{0 & \mu_1}} & R/x \ar[d]^{\smat{a \mu_x \\ \mu_x}} \ar[r]^-{\smat{0 \\ \mu_x}} & R/x \op R/x^2 \ar@{=}[d] \\
R/x^3 \op R/x^2 \ar[r]_-{\smat{\mu_1 & 0}} & R/x \ar[r]_-{\smat{\mu_x \\ 0}} & R/x^2 \op R/x^2 \ar[r]_-{\smat{\mu_1 & 0 \\ 0 & 1}} & R/x \op R/x^2, \\
}
\end{equation}
a general fill-in has the stated form, with $a \in \F_2$. By \cref{lem:IntoProduct}, the map \eqref{eq:F2C4good} is good if and only if its restriction and projection
\[
\cxymatrix{
R/x^2 \ar[d]_{0} \ar[r]^-{\mu_x} & R/x^3 \ar[d]^{0} \ar[r]^-{\mu_1} & R/x \ar[d]^{a \mu_x} \ar[r]^-{\mu_x} & R/x^2 \ar[d]^{0} \\
R/x^3 \ar[r]_-{\mu_1} & R/x \ar[r]_-{\mu_x} & R/x^2 \ar[r]_-{\mu_1} & R/x \\
}
\]
is good. By \cite{Neeman91}*{\S 1, Case~2}, this holds if and only if $a \mu_x$ is a lightning flash, which holds if and only if $a=0$ holds. 
Starting with the non-good fill-in $\smat{\mu_x \\ \mu_x}$ %
and using the 
shear automorphism
\[
\smat{1 & 1 \\ 0 & 1} \colon R/x^2 \op R/x^2 \ral{\cong} R/x^2 \op R/x^2
\]
in the bottom row, we obtain a map of triangles
\begin{equation}\label{eq:F2C4notGood}
\cxymatrix{
R/x^3 \op R/x^2 \ar@{=}[d] \ar[r]^-{\smat{1 & 0 \\ 0 & \mu_x}} & R/x^3 \op R/x^3 \ar[d]^{\smat{\mu_1 & 0}} \ar[r]^-{\smat{0 & \mu_1}} & R/x \ar[d]^{\smat{0 \\ \mu_x}} \ar[r]^-{\smat{0 \\ \mu_x}} & R/x \op R/x^2 \ar@{=}[d] \\
R/x^3 \op R/x^2 \ar[r]_-{\smat{\mu_1 & 0}} & R/x \ar[r]_-{\smat{\mu_x \\ 0}} & R/x^2 \op R/x^2 \ar[r]_-{\smat{\mu_1 & \mu_1 \\ 0 & 1}} & R/x \op R/x^2 \\
}
\end{equation}
which is \emph{not} good and has the same (homotopy cartesian) middle square as the good map \eqref{eq:F2C4good} in the case $a=0$. 
Alternately, one can check the (non-)goodness claims using \cref{prop:HoCartFillin}~\eqref{item:1ghGood}. %

A straightforward calculation shows that the map~\eqref{eq:F2C4notGood} does not extend to an octahedron, so it is \emph{not} Verdier good, by \cref{pr:1ghVerdier}. 
Indeed, \eqref{eq:F2C4notGood} does extend to a $4 \x 4$ diagram
\[
\cxymatrix{
R/x^3 \op R/x^2 \ar@{=}[d] \ar[r]^-{\smat{1 & 0 \\ 0 & \mu_x}} & R/x^3 \op R/x^3 \ar[d]^{\smat{\mu_1 & 0}} \ar[r]^-{\smat{0 & \mu_1}} & R/x \ar[d]^{\smat{0 \\ \mu_x}} \ar[r]^-{\smat{0 \\ \mu_x}} & R/x \op R/x^2 \ar@{=}[d] \\
R/x^3 \op R/x^2 \ar[r]_-{\smat{\mu_1 & 0}} & R/x \ar[d]^{\smat{\mu_x \\ 0}} \ar[r]_-{\smat{\mu_x \\ 0}} & R/x^2 \op R/x^2 \ar[d]^{\smat{1 & 0 \\ 0 & \mu_1}} \ar[r]_-{\smat{\mu_1 & \mu_1 \\ 0 & 1}} & R/x \op R/x^2 \\
& R/x^2 \op R/x \ar[d]^{\smat{\mu_1 & 0 \\ c \mu_1 & 1}} \ar@{=}[r] & R/x^2 \op R/x \ar[d]^{\smat{0 & \mu_{x^2}}} & \\
& R/x \op R/x \ar[r]_-{\smat{0 & \mu_{x^2}}} & R/x^3 & \\
}
\]
with $c \in \F_2$,
and, up to isomorphism, every extension is of this form.
But this diagram does not satisfy the
additional equation required of an octahedron:
\begin{align*}
&\left( \Si \smat{1 & 0 \\ 0 & \mu_x} \right) \smat{\mu_1 & \mu_1 \\ 0 & 1} 
= \smat{1 & 0 \\ 0 & \mu_1} \smat{\mu_1 & \mu_1 \\ 0 & 1} 
= \smat{\mu_1 & \mu_1 \\ 0 & \mu_1} \\
\neq &\smat{\mu_1 & 0 \\ c \mu_1 & 1} \smat{1 & 0 \\ 0 & \mu_1} 
= \smat{\mu_1 & 0 \\ c \mu_1 & \mu_1} \colon R/x^2 \op R/x^2 \to R/x \op R/x.
\end{align*}
So \eqref{eq:F2C4notGood} is middling good, but not Verdier good.
This also provides an alternate proof that the map~\eqref{eq:F2C4notGood} is not good, via \cref{cor:1ghVerdier}.
\end{ex}

\begin{ex}\label{ex:HoCartNotGood}
Consider the derived category of the integers $\DZ$, with the notation as in \cref{ex:naive}. 
Let $n \geq 2$ be an integer.  
Using \cref{prop:HoCartFillin}~\eqref{item:1ghGood}, one can show that the map of triangles
\[
\cxymatrix{
\Z[0] \ar@{=}[d] \ar[r]^-{2} & \Z[0] \ar[d]^{q_4} \ar[r]^-{q_2} & \Z/2[0] \ar[d]^{\smat{1 \\ 0}} \ar[r]^-{\ep_2} & \Z[1] \ar@{=}[d] \\
\Z[0] \ar[r]_-{2q_4} & \Z/4[0] \ar[r]_-{\smat{q \\ -\ep_4}} & \Z/2[0] \op \Z[1] \ar[r]_-{\smat{\ep_2 & 2}} & \Z[1]\\
}
\]
is good. In fact, $\smat{1 \\ 0}$ is the unique good fill-in here.
However, the map of triangles with the same %
middle square 
\begin{equation}\label{eq:DZnotGood}
\cxymatrix{
\Z[0] \ar@{=}[d] \ar[r]^-{2} & \Z[0] \ar[d]^{q_4} \ar[r]^-{q_2} & \Z/2[0] \ar[d]^{\smat{1 \\ 0}} \ar[r]^-{\ep_2} & \Z[1] \ar@{=}[d] \\
\Z[0] \ar[r]_-{2q_4} & \Z/4[0] \ar[r]_-{\smat{q \\ -\ep_4}} & \Z/2[0] \op \Z[1] \ar[r]_-{\smat{\ep_2 & -2}} & \Z[1] \\
}
\end{equation}
is \emph{not} good. In fact, the unique good fill-in is $\smat{1 \\ \ep_2}$. 

As in \cref{ex:F2C4}, the map~\eqref{eq:DZnotGood} extends to a $4 \x 4$ diagram, but \emph{not} to an octahedron, so it is middling good but not Verdier good.
\end{ex}

\begin{ex}\label{ex:F2C4composite}
Let us return to the stable module category $\StMod(R)$ %
from \cref{ex:F2C4}. Consider the diagram with exact rows
\[
\resizebox{\textwidth}{!}{
$\xymatrix @C+0.5pc {
R/x^3 \op R/x^2 \op R/x^2 \ar@{=}[d] \ar[r]^-{\smat{1 & 0 & 0 \\ 0 & 1 & 0 \\ 0 & 0 & \mu_x}} & R/x^3 \op R/x^2 \op R/x^3 \ar[d]^{\smat{\mu_1 & 0 & 0 \\ 0 & 0 & 1}} \ar[r]^-{\smat{0 & 0 & \mu_1}} & R/x \ar[d]^{\smat{0 \\ 0 \\ 1}} \ar[r]^-{\smat{0 \\ 0 \\ \mu_x}} & R/x \op R/x^2 \op R/x^2 \ar@{=}[d] \\
R/x^3 \op R/x^2 \op R/x^2 \ar@{=}[d] \ar[r]^-{\smat{\mu_1 & 0 & 0 \\ 0 & 0 & \mu_x}} & R/x \op R/x^3 \ar[d]^{\smat{1 & 0}} \ar[r]^-{\smat{\mu_x & 0 \\ 0 & 0 \\ 0 & \mu_1}} & R/x^2 \op R/x^2 \op R/x \ar[d]^{\smat{1 & 0 & 0 \\ 0 & 1 & 0 \\ 0 & 0 & \mu_x}} \ar[r]^-{\smat{\mu_1 & 0 & 0 \\ 0 & 1 & 0 \\ 0 & 0 & \mu_x}} & R/x \op R/x^2 \op R/x^2 \ar@{=}[d] \\
R/x^3 \op R/x^2 \op R/x^2 \ar[r]_-{\smat{\mu_1 & 0 & 0}} & R/x \ar[r]_-{\smat{\mu_x \\ 0 \\ 0}} & R/x^2 \op R/x^2 \op R/x^2 \ar[r]_-{\smat{\mu_1 & 0 & \mu_1 \\ 0 & 1 & 0 \\ 0 & 0 & 1}} & R/x \op R/x^2 \op R/x^2. \\
}$
}
\]
By \cref{lem:IntoProduct} and \cref{pr:IncUpperTriangular}~\eqref{item:Inc}, the top part is a good map of triangles. 
By \cref{lem:IntoProduct} and \cref{pr:IncUpperTriangular}~\eqref{item:Proj}, %
the bottom part is a good map of triangles. 
By the same argument as in \cref{ex:F2C4}, 
the composite map of triangles from the top row to the bottom row is \emph{not} good. 
\end{ex}

\section{Maps of triangles with one or two zero components}\label{se:zero}

In this section we study maps of triangles with at most two non-zero components
and show that the properties of being good, Verdier good,
and nullhomotopic agree for such maps in most cases.
This leads to a lifting condition expressed in terms of Verdier good fill-ins.

\begin{lem}\label{lem:0gh}
A map $(0,g,h)$ of triangles is nullhomotopic if and only if it is nullhomotopic
via a chain homotopy $(0,G,0)$.  Rotating the triangles gives similar statements
for any map with at most two non-zero components.
\end{lem}

\begin{proof}
Assume that the map $(0,g,h)$ is nullhomotopic via a chain homotopy $(F,G,H)$, as illustrated in the diagram
\[
\cxymatrix{
X \ar[d]_{0} \ar[r]^-{u} & Y \ar[dl]_{F} \ar[d]^{g} \ar[r]^-{v} & Z \ar[dl]_{G} \ar[d]^{h} \ar[r]^-{w} & \Sigma X \ar[dl]_{H} \ar[d]^{0} \\
X' \ar[r]_-{u'} & Y' \ar[r]_-{v'} & Z' \ar[r]_-{w'} & \Sigma X'. \\
}
\]
In particular, the map $0 \colon \Si X \to \Si X'$ satisfies $0 = w'H + \Si (Fu)$. Hence, the following solid diagram commutes:
\[
\xymatrix @C+0.4pc {
X \ar[d]_{\Si^{-1} H} \ar[r]^-{u} & Y \ar[d]^{F} \ar[r]^-{v} & Z \ar@{-->}[d]^{\tild{G}} \ar[r]^-{w} & \Sigma X \ar[d]^{H} \\
\Si^{-1} Z' \ar[r]_-{-\Si^{-1} w'} & X' \ar[r]_-{u'} & Y' \ar[r]_-{v'} & Z'. \\
}
\]
Let $\tild{G} \colon Z \to Y'$ be a fill-in. Then $(0, G+\tild{G}, 0)$ is a nullhomotopy of $(0,g,h)$:
\begin{align*}
(G+\tild{G})v + u'0 &= Gv + u'F = g \\
0w + v'(G+\tild{G}) &= v'G + Hw = h. \qedhere
\end{align*}
\end{proof}

\begin{cor}\label{cor:Nullhomot1}
A map of triangles with only one non-zero component, i.e., of the form $(f,0,0)$ or $(0,g,0)$ or $(0,0,h)$, is nullhomotopic if and only if that component is a lightning flash. \qed
\end{cor}

\begin{rem}
By considering the difference of two maps of triangles, \cref{lem:0gh} (resp.\ \cref{cor:Nullhomot1}) characterizes the chain homotopies that modify only two components (resp.\ only one component).
\end{rem}

\begin{prop}\label{pr:0ghVerdierGood}
For a map of triangles with at most two non-zero components, the following are equivalent.
\begin{enumerate}
 \item \label{item:Nullhomotopic} The map is nullhomotopic.
 \item \label{item:Good} The map is good.
\end{enumerate}
If the map has the form $(0,g,h)$ or $(f,0,h)$, then those conditions are further equivalent to:
\begin{enumerate}
\setcounter{enumi}{2}
 \item \label{item:VerdierGood} The map is Verdier good.
\end{enumerate}
\end{prop}

We don't know whether being middling good is equivalent to the above conditions.
We consider maps of the form $(f,g,0)$ in \cref{pr:fg0Good}.

\begin{proof}
\Eqref{item:Nullhomotopic}~$\iff$~\Eqref{item:Good}. Since goodness and the property of being nullhomotopic are invariant under rotation, it suffices to treat one of the three cases, say, $(f,g,0)$. By \cref{th:obstruction}, if the map $(f,g,0)$ is good, then it is nullhomotopic via a chain homotopy $(F,0,0)$. Conversely, any nullhomotopic map is good. 

Since Verdier goodness is not known to be rotation invariant, we will %
prove the equivalence \Eqref{item:VerdierGood}~$\iff$~\Eqref{item:Nullhomotopic} 
in separate cases.

\textbf{Case $(0,g,h)$.} An octahedron for the composite $X \ral{u} Y \ral{g} Y'$ is of the form: 
\begin{equation}\label{eq:OctahedronUG0}
\cxymatrix{
X \ar@{=}[d] \ar[r]^-{u} & Y \ar[d]^{g} \ar[r]^-{v} & Z \ar@{-->}[d]^{\al_1 = \smat{a_1 \\ w\;}} \ar[r]^-{w} & \Si X \ar@{=}[d] \\
X \ar[r]^-{0} & Y' \ar[d]^{g'} \ar[r]^-{\inc} & Y' \op \Si X \ar@{-->}[d]^{\be_1 = \smat{g' & b_1}} \ar[r]^-{\proj} & \Si X \\
& Y'' \ar[d]^{g''} \ar@{=}[r] & Y'' \ar@{-->}[d]^{\ga_1 = (\Si v)g''} & \\
& \Si Y \ar[r]_-{\Si v} & \Si Z, & \\
}
\end{equation}
with $a_1 v = g$ %
and $g'' b_1 = \Si u$. %
An octahedron for the composite $X \ral{0} X' \ral{u'} Y'$ is of the form: 
\begin{equation}\label{eq:Octahedron0U0}
\cxymatrix{
X \ar@{=}[d] \ar[r]^-{0} & X' \ar[d]^{u'} \ar[r]^-{\inc} & X' \op \Si X \ar@{-->}[d]^{\al_2 = \smat{u' & a_2 \\ 0\; & 1}} \ar[r]^-{\proj} & \Si X \ar@{=}[d] \\
X \ar[r]^-{0} & Y' \ar[d]^{v'} \ar[r]^-{\inc} & Y' \op \Si X \ar@{-->}[d]^{\be_2 = \smat{v' & b_2}} \ar[r]^-{\proj} & \Si X \\
& Z' \ar[d]^{w'} \ar@{=}[r] & Z' \ar@{-->}[d]^{\ga_2 = \smat{w' \\ 0\;}} & \\
& \Si X' \ar[r]_-{\inc} & \Si X' \op \Si^2 X, & \\
}
\end{equation}
with $w' b_2 = 0$. %
Verdier good fill-ins are the maps
\[
h = \be_2 \circ \al_1 = v' a_1 + b_2 w.
\]
For such a fill-in $h$, the map of triangles $(0,g,h)$ is nullhomotopic via the chain homotopy $(0, a_1, b_2)$, as illustrated in the diagram
\[
\xymatrix{
X \ar[d]_{0} \ar[r]^-{u} & Y \ar[dl]_{0} \ar[d]^{g} \ar[r]^-{v} & Z \ar[dl]_{a_1} \ar[d]^{h} \ar[r]^-{w} & \Sigma X \ar[d]^{0} \ar[dl]_{b_2} \\
X' \ar[r]_-{u'} & Y' \ar[r]_-{v'} & Z' \ar[r]_-{w'} & \Sigma X'. \\
}
\]

Conversely, assume that $(0,g,h)$ is nullhomotopic. By \cref{lem:0gh}, we may assume that the nullhomotopy is of the form $(0,a_1,0)$, 
so that in particular %
$h = v' a_1$ holds. 
We will show that the given $a_1$ and %
$b_2 = 0$ 
arise as valid choices in the octahedra above.
By \cref{prop:HoCartFillin}, there exists a fill-in in the diagram
\[
\cxymatrix{
Y \ar@{=}[d] \ar[r]^-{v} & Z \ar[d]_{a_1} \ar[r]^-{w} & \Si X \ar@{-->}[d]^{- b_1} \ar[r]^-{- \Si u} & \Si Y \ar@{=}[d] \\
Y \ar[r]_-{g} & Y' \ar[r]_-{g'} & Y'' \ar[r]_-{g''} & \Si Y \\
}
\]
making the middle square homotopy cartesian with differential $\del = (\Si v) g''$, i.e., making the following triangle exact:
\[
\xymatrix @C+0.3pc {
Z \ar[r]^-{\smat{a_1 \\ w\;}} & Y' \op \Si X \ar[r]^-{\smat{g' & b_1}} & Y'' \ar[r]^{(\Si v) g''} & \Si Z. \\
}
\]
This is the third column of 
\eqref{eq:OctahedronUG0}. Moreover $b_1$ satisfies $g'' b_1 = \Si u$, as required. 

For~\eqref{eq:Octahedron0U0}, picking $a_2=0$ and $b_2=0$ yields a valid choice. 
It follows that %
$h = v' a_1$ 
is a Verdier good fill-in.

\textbf{Case $(f,0,h)$.} The argument is the same as the previous case, using a good fibre fill-in in \cref{prop:HoCartFillin}.
\end{proof}

Now we can refine \cref{th:obstruction}.

\begin{cor}[Lifting Criterion]\label{cor:lifting}
Given a solid arrow commutative square
\[
  \cxymatrix{
    X \ar[r]^u \ar[d]_f & Y \ar[d]^g \ar@{-->}[dl]_k \\
    X' \ar[r]_{u'} & Y' , \\
  }
\]
choose extensions of the rows to cofibre sequences as in~\eqref{eq:fill-ins}. The following are equivalent.
\begin{enumerate}
\item There exists a lift $k \colon Y \to X'$.
\item The map $0 \colon Z \to Z'$ is a fill-in and the map $(f,g,0)$ of triangles is good.
\item The map $0 \colon Z \to Z'$ is a fill-in and the map $(0,f,g)$ of rotated triangles is Verdier good. \qed
\end{enumerate}
\end{cor}

We need to rotate the triangles in the last case because \cref{pr:0ghVerdierGood}
does not include the case $(f,g,0)$.

As with \cref{th:obstruction}, one can use \cref{cor:lifting} to define an
obstruction class consisting of the fill-in maps $Z \to Z'$ such that the
rotated map $(\Sigma^{-1} h, f, g)$ is Verdier good.
A lift exists if and only if this obstruction class contains zero.

Next, we tackle the remaining case $(f,g,0)$, obtaining a weaker result in \cref{pr:fg0Good}.

\begin{defn}\label{def:WeaklyFunctorial}
A triangulated category $\cT$ admits \Def{weakly functorial octahedra} if for every diagram
\begin{equation}\label{eq:MapBases}
\cxymatrix{
X_1 \ar[d]_{\phy_1} \ar[r]^-{f_1} & X_2 \ar[d]^{\phy_2} \ar[r]^-{f_2} & X_3 \ar[d]^{\phy_3} \\
Y_1 \ar[r]_-{g_1}  & Y_2 \ar[r]_-{g_2} & Y_3 \\
}
\end{equation}
in $\cT$, there exist octahedra $\mathcal{X}$ and $\mathcal{Y}$ based on the top and bottom rows of \eqref{eq:MapBases} respectively, and a map of octahedra $\phy \colon \mathcal{X} \to \mathcal{Y}$ extending the given map between their bases $(\phy_1, \phy_2, \phy_3)$.
\end{defn}

\begin{ex}\label{ex:n-angulated}
If a triangulated category $\cT$ admits a $3$-pretriangulated enhancement in the sense of Maltsiniotis \cite{Maltsiniotis06}, then it admits weakly functorial octahedra. Indeed, pick \emph{distinguished} octahedra (i.e., $3$-triangles) $\mathcal{X}$ and $\mathcal{Y}$ based on $X_1 \to X_2 \to X_3$ and $Y_1 \to Y_2 \to Y_3$ respectively. Then any map $(\phy_1, \phy_2, \phy_3)$ between their bases as in~\eqref{eq:MapBases} extends to a map of octahedra $\phy \colon \mathcal{X} \to \mathcal{Y}$. This is stronger than the condition in \cref{def:WeaklyFunctorial}, where the two octahedra may depend on the given map between bases. %

The homotopy category of a stable model category (or more generally of a complete and cocomplete stable $\infty$-category) admits a canonical $n$-triangulation for every $n$; see \cite{Maltsiniotis06} and \cite{GrothS16}*{Theorem~13.6, Examples~13.8}. 
See also related work on higher triangulations in \cite{Kunzer07}. %
\end{ex}

\begin{prop}\label{pr:fg0Good}
Assume that the triangulated category $\cat{T}$ admits weakly functorial octahedra. If a map of triangles $(f,g,0)$ is good, then it is Verdier good.
\end{prop}

\begin{proof}
Consider the diagram:
\begin{equation}\label{eq:MapBasesNull}
\cxymatrix{
X \ar@{=}[d] \ar[r]^-{u} & Y \ar[d]^F \ar[r]^-{g} & Y' \ar@{=}[d] \\
X \ar[r]_-{f} & X' \ar[r]_-{u'} & Y'. \\
}
\end{equation}
By assumption, we may pick octahedra based on $X \ral{u} Y \ral{g} Y'$ and $X \ral{f} X' \ral{u'} Y'$ that admit a map of octahedra $\phy$ extending the map between their bases~\eqref{eq:MapBasesNull}. 
Keeping the notation as in \cref{def:VerdierGood}, denote the new components of the natural transformation $\phy$ by $\phy_Z \colon Z \to X''$, $\phy_A \colon A \to A$, and $\phy_{Y''} \colon Y'' \to Z'$. Due to the map of triangles
\[
\xymatrix{
X \ar@{=}[d] \ar[r]^-{gu} & Y' \ar@{=}[d] \ar[r]^-{\tild{v}} & A \ar[d]^{\phy_A} \ar[r]^-{\tild{w}} & \Si X \ar@{=}[d] \\
X \ar[r]_-{u'f} & Y' \ar[r]_-{\tild{v}} & A \ar[r]_-{\tild{w}} & \Si X, \\
}
\]
the map $\phy_A$ is an isomorphism. Moreover, the maps $Z \ral{\phy_A \al_1} A \ral{\be_1 \phy_A^{-1}} Y''$ remain valid choices for the first octahedron. These choices yield a Verdier good fill-in
\[
h = \be_2 (\phy_A \al_1) = \be_2 \al_2 \phy_Z = 0 \phy_Z = 0,
\]
where the second equality uses the naturality of $\phy$.
\end{proof}

We now study a particularly simple family of maps of triangles.

\begin{prop}\label{pr:ZeroVerdierGood}
For a map of triangles of the form
\begin{equation*}
\cxymatrix{
X \ar[d]_{0} \ar[r]^-{u} & Y \ar[d]^{0} \ar[r]^-{v} & Z \ar[d]^{h} \ar[r]^-{w} & \Sigma X \ar[d]^{0} \\
X' \ar[r]_-{u'} & Y' \ar[r]_-{v'} & Z' \ar[r]_-{w'} & \Sigma X' , \\
}
\end{equation*}
the following are equivalent:
\begin{enumerate}
\item \label{item:hGood} The map $(0,0,h)$ is good.
\item \label{item:hVerdier} The map $(0,0,h)$ is Verdier good.
\item \label{item:Middling} The map $(0,0,h)$ is middling good.
\item \label{item:Toda} The Toda bracket $\lan w', h, v \ran \subseteq \cat{T}(\Si Y, \Si X')$ contains zero.
\item \label{item:Lightning} The map $h$ is a lightning flash.
\end{enumerate}
Moreover, analogous statements hold for maps of the form $(f,0,0)$ and $(0,g,0)$.
Since the other conditions are rotation invariant, it follows that Verdier goodness
is rotation invariant for maps with at most one non-zero component.
\end{prop}

\begin{proof}
The equivalence of~\Eqref{item:hGood}, \Eqref{item:hVerdier} and~\Eqref{item:Lightning}
follows from \cref{cor:Nullhomot1} and \cref{pr:0ghVerdierGood}.
(Neeman also showed that \Eqref{item:hGood} is equivalent to~\Eqref{item:Lightning}
in \cite{Neeman91}*{\S 1, Case~2}.)
That \Eqref{item:hVerdier} implies \Eqref{item:Middling} is exactly Verdier's argument, which we
summarized in \cref{lem:VerdierHaynes}.

We will show the implications~\Eqref{item:Middling} $\implies$ \Eqref{item:Toda} $\implies$ \Eqref{item:Lightning}.

To show that~\Eqref{item:Middling} implies~\Eqref{item:Toda}, assume that $(0,0,h)$ is middling good.
A $4 \x 4$ diagram extending the map $(0,0,h)$ may be assumed to be of the form
\[
\xymatrix @C+1.2pc @R+0.2pc {
X \ar[d]_{0} \ar[r]^-{u} & Y \ar[d]^{0} \ar[r]^-{v} & Z \ar[d]^h \ar[r]^-{w} & \Si X \ar[d]^{0} \\
X' \ar[d]_{\inc} \ar[r]^-{u'} & Y' \ar[d]^{\inc} \ar[r]^-{v'} & Z' \ar[d]^(.4){h'} \ar[r]^{w'} & \Si X' \ar[d]^{\inc} \\
X' \op \Si X \ar[d]_{\proj} \ar@{-->}[r]^-{\smat{u' & a \\ 0\; & \Si u}} & Y' \op \Si Y \ar[d]^{\proj} \ar@{-->}[r]^-{\smat{h'v' & b\;}} & Z'' \anti \ar[d]^{h''} \ar@{-->}[r]^-{\smat{c\;\;\;\; \\ -(\Si w) h''}} & \Si X' \op \Si^2 X \ar[d]^{\proj} \hspace*{-8pt}\\
\Si X \ar[r]^-{\Si u} & \Si Y \ar[r]^-{\Si v} & \Si Z \ar[r]^-{\Si w} & \Si^2 X \\
}
\]
with $h'' b = \Si v$ and $c h' = w'$. 
These equations exhibit the composite $c(-b) \in \lan w', h, v \ran$ as being in the fibre-cofibre Toda bracket \cite{ChristensenFrankland17}*{Definition~3.1}.
Exactness of the third row implies that $c b = 0$. 
Hence, $0 \in \lan w', h, v \ran$.

To show that~\Eqref{item:Toda} implies~\Eqref{item:Lightning}, assume that $\lan w', h, v \ran$ contains zero.
Using the iterated cofibre Toda bracket \cite{ChristensenFrankland17}*{Definition~3.1},
this means that there exists a map $\phy$ making the diagram
\[
\cxymatrix{
Y \ar@{=}[d] \ar[r]^-{v} & Z \ar@{=}[d] \ar[r]^-{w} & \Si X \ar[d]^{\phy} \ar[r]^-{-\Si u} & \Sigma Y \ar[d]^{0} \\
Y \ar[r]_-{v} & Z \ar[r]_-{h} & Z' \ar[r]_-{w'} & \Sigma X'\\
}
\]
commute.
The condition $w' \phy = 0$ implies that $\phy$ lifts as $\phy = v' \tild{\phy}$ for some $\tild{\phy} \colon \Si X \to Y'$.  So we have $h = \phy w = v' \tild{\phy} w$, which is a lightning flash.
\end{proof}

It follows that we can characterize when an arbitrary 3-fold Toda bracket
contains zero: for any composable maps $Y \ral{v} Z \ral{h} Z' \ral{w'} W$,
the Toda bracket $\lan w', h, v \ran$ contains zero if and only if
$h = v' \te w$ for some $\te$, where $v'$ is the fibre of $w'$ and $w$ is the cofibre of $v$.

\section{Maps between rotations of a triangle}\label{se:rotations}

In this final section, we study a special situation in which we again find
that good morphisms and Verdier good morphisms agree, and use this to
present some interesting examples.

\begin{prop}\label{pr:FillInRotated}
Consider a map of triangles of the form
\[
\cxymatrix{
X \ar[d]_{u} \ar[r]^-{u} & Y \ar[d]^{v} \ar[r]^-{v} & Z \ar[d]^h \ar[r]^-{w} & \Sigma X \ar[d]^{\Si u} \\
Y \ar[r]_-{v} & Z \ar[r]_-{w} & \Si X \ar[r]_-{-\Si u} & \Sigma Y, \\
}
\]
where the second row is the displayed rotation of the first row, and the
first two vertical maps are taken from the triangle.
Then the following are equivalent:
\begin{enumerate}
  \item The map $(u,v,h)$ is good.
  \item The map $(u,v,h)$ is Verdier good.
  \item The Toda bracket $\lan -\Si u, h, v \ran \subseteq \cat{T}(\Si Y, \Si Y)$ contains zero.
  \item The map $h$ is a lightning flash.
\end{enumerate}
When $h = w$, these conditions are also equivalent to the given triangle being contractible.
\end{prop}

\begin{proof}
A chain homotopy consisting of the identity map $1_Y$ and zero elsewhere shows
that $(u,v,h)$ is chain homotopic to $(0,0,h)$.
Since goodness is invariant under chain homotopy (\cref{rem:good-homotopic})
condition~(1) is equivalent to $(0,0,h)$ being good.
By \cref{pr:ZeroVerdierGood}, this is equivalent to conditions~(3) and~(4).
In the case that $h = w$, Neeman showed in~\cite{Neeman91}*{\S 1, Case~2}
that this is also equivalent to the triangle being contractible.
We will show that (2) is equivalent to (4).

We begin by assuming that $(u,v,h)$ is Verdier good.
Every octahedron for the composite \mbox{$X \ral{u} Y \ral{v} Z$} is isomorphic to one of the form:
\begin{equation}\label{eq:OctahedronUV}
\cxymatrix{
X \ar@{=}[d] \ar[r]^-{u} & Y \ar[d]^{v} \ar[r]^-{v} & Z \ar@{-->}[d]^{\al = \smat{a \\ w}} \ar[r]^-{w} & \Si X \ar@{=}[d] \\
X \ar[r]^-{0} & Z \ar[d]^{w} \ar[r]^-{\inc} & Z \op \Si X \ar@{-->}[d]^{\be = \smat{w & b}} \ar[r]^-{\proj} & \Si X \\
& \Si X \ar[d]^{-\Si u} \ar@{=}[r] & \Si X \ar@{-->}[d]^{\ga=0} & \\
& \Si Y \ar[r]_-{\Si v} & \Si Z. & \\
}
\end{equation}
The top-middle square gives the equation $av = v$, that is, $(a-1)v = 0$.
Therefore, $a = 1 + \phy w$ for some $\phy : \Si X \to Z$.
The staircase equation from $Z \op \Si X$ to $\Si Y$ is
$(\Si u) \circ \proj = (-\Si u) \circ \begin{bmatrix} w & b \end{bmatrix}$,
or equivalently, $(\Si u) (b+1) = 0$.
Therefore, $b = -1 + w \psi$ for some $\psi : \Si X \to Z$.

Making such choices for two octahedra, every Verdier good fill-in will be of the form
\begin{align*}
h = \be_2 \al_1 &= \begin{bmatrix} w & b_2 \end{bmatrix} \begin{bmatrix} a_1 \\ w \end{bmatrix} \\
&= w a_1 + b_2 w \\
&= w (1 + \phy_1 w) + (-1 + w \psi_2) w \\
&= w (\phy_1 + \psi_2) w ,
\end{align*}
which is a lightning flash, proving (4).

To prove that (4) implies (2), assume that $h$ is a lightning flash.
By \cref{lem:AddLightning}, it is enough to show that $(u,v,0)$ is Verdier good.
One can check that choosing $a = 1$ and $b = -1$ in~\eqref{eq:OctahedronUV} gives an
octahedron.  Making these choices for both octahedra, the Verdier fill-in is
\[
  \be_2 \al_1
  = \begin{bmatrix} w & -1 \end{bmatrix} \begin{bmatrix} 1 \\ w \end{bmatrix}
  = 0 .
\]
So $(u,v,0)$ is Verdier good.
\end{proof}

\begin{ex}\label{ex:MiddlingNotHoInvt}
As mentioned in the proof of \cref{pr:FillInRotated},
Neeman shows in \cite{Neeman91}*{\S 1, Case~2} that the map of triangles
\[
\cxymatrix{
X \ar[d]_{0} \ar[r]^-{u} & Y \ar[d]^{0} \ar[r]^-{v} & Z \ar[d]^w \ar[r]^-{w} & \Sigma X \ar[d]^{0} \\
Y \ar[r]_-{v} & Z \ar[r]_-{w} & \Si X \ar[r]_-{-\Si u} & \Sigma Y \\
}
\]
is good if and only if the given triangle is contractible.
By \cref{pr:ZeroVerdierGood}, this is equivalent to the map $(0,0,w)$ being middling good.
However, as we observed in the proof of \cref{pr:FillInRotated},
this map of triangles is chain homotopic to the map
\[
\cxymatrix{
X \ar[d]_{u} \ar[r]^-{u} & Y \ar[d]^{v} \ar[r]^-{v} & Z \ar[d]^w \ar[r]^-{w} & \Sigma X \ar[d]^{\Si u} \\
Y \ar[r]_-{v} & Z \ar[r]_-{w} & \Si X \ar[r]_-{-\Si u} & \Sigma Y. \\
}
\]
This map $(u,v,w)$ is always middling good, but is good if and only if it is
Verdier good if and only if the given triangle is contractible, by 
\cref{pr:FillInRotated}.
\end{ex}

As a consequence, we deduce the following result.

\begin{cor}\label{cor:examples}
\begin{enumerate}
\item There exist maps of triangles which are not middling good.
\item There exist maps of triangles which are middling good but neither
      good nor Verdier good.
\item Middling goodness is not invariant under chain homotopy.
\end{enumerate}
\end{cor}

\begin{proof}
Consider a non-contractible triangle in \cref{ex:MiddlingNotHoInvt}.
(We will give some non-contractible triangles in \cref{ex:DZ,ex:Spectra}.)
Then $(0,0,w)$ is not middling good, proving (1).
And $(u,v,w)$ is middling good, but is neither good nor Verdier good, proving (2).
These two maps are chain homotopic, proving~(3).
\end{proof}

\begin{ex}\label{ex:DZ}
In the derived category of the integers $\DZ$, consider the rotation of the map of triangles in \cref{ex:naive}: 
\[
\xymatrix{
\Z[0] \ar[d]_{0} \ar[r]^-{n} & \Z[0] \ar[d]^{0} \ar[r]^-{q_n} & \Z/n[0] \ar[d]^{\ep_n} \ar[r]^-{\ep_n} & \Z[1] \ar[d]^{0} \\
\Z[0] \ar[r]_-{q_n} & \Z/n[0] \ar[r]_-{\ep_n} & \Z[1] \ar[r]_-{-n} & \Z[1]. \\
}
\]
All maps $\Z[1] \to \Z/n[0]$ are zero for degree reasons, so the fill-in $\ep_n$ is not a lightning flash. 
It follows that this triangle is not contractible, that $(0,0,\ep_n)$ is not middling good, etc.
One can also show directly that $(0,0,\ep_n)$ does not extend to a $4 \times 4$ diagram of triangles.
\end{ex}

\begin{ex}\label{ex:Spectra}
In the stable homotopy category, consider the mod $n$ Moore spectrum $M(n)$ for $n \geq 2$,
sitting in an exact triangle
\[
\xymatrix{
S^0 \ar[r]^-{n} & S^0 \ar[r]^-{q} & M(n) \ar[r]^-{\de} & S^1. \\
}
\]
By the long exact sequence of homotopy groups, the connecting map $\de$ induces
\[
0 = \pi_1(\de) \colon \pi_1 M(n) \to \pi_1 S^1.
\]
It follows that $\delta$ is not a lightning flash when regarded as part of a
map $(0,0,\delta)$ of triangles.
Again, this implies that this triangle is not contractible, etc.
\end{ex}

\vspace*{10pt}

\end{document}